\documentclass[a4paper,leqno,10pt]{amsart}

\usepackage{epsf,graphicx,epsfig}
\usepackage{amscd}
\usepackage{amsmath,latexsym,amssymb,amsthm}
\usepackage[nospace,noadjust]{cite}
\usepackage{textcomp}
\usepackage{setspace,cite}
\usepackage{lscape,fancyhdr,fancybox}
\usepackage{stmaryrd}
\usepackage[all,cmtip]{xy}
\RequirePackage{tikz-cd}
\RequirePackage{amssymb}
\usetikzlibrary{calc}
\usetikzlibrary{decorations.pathmorphing}

\tikzset{curve/.style={settings={#1},to path={(\tikztostart)
			.. controls ($(\tikztostart)!\pv{pos}!(\tikztotarget)!\pv{height}!270:(\tikztotarget)$)
			and ($(\tikztostart)!1-\pv{pos}!(\tikztotarget)!\pv{height}!270:(\tikztotarget)$)
			.. (\tikztotarget)\tikztonodes}},
	settings/.code={\tikzset{quiver/.cd,#1}
		\def\pv##1{\pgfkeysvalueof{/tikz/quiver/##1}}},
	quiver/.cd,pos/.initial=0.35,height/.initial=0}

\tikzset{tail reversed/.code={\pgfsetarrowsstart{tikzcd to}}}
\tikzset{2tail/.code={\pgfsetarrowsstart{Implies[reversed]}}}
\tikzset{2tail reversed/.code={\pgfsetarrowsstart{Implies}}}
\tikzset{no body/.style={/tikz/dash pattern=on 0 off 1mm}}

\usepackage{graphicx}

\usepackage{color}
\usepackage{url}
\usepackage{enumerate}
\usepackage[mathscr]{euscript}

\newcounter{cnt1}
\newcounter{cnt2}
\newcounter{cnt3}
\newcommand{\blr}{\begin{list}{$($\roman{cnt1}$)$} {\usecounter{cnt1}
			\setlength{\topsep}{0pt} \setlength{\itemsep}{0pt}}}
	\newcommand{\bla}{\begin{list}{$($\alph{cnt2}$)$} {\usecounter{cnt2}
				\setlength{\topsep}{0pt} \setlength{\itemsep}{0pt}}}
		\newcommand{\bln}{\begin{list}{$($\arabic{cnt3}$)$} {\usecounter{cnt3}
					\setlength{\topsep}{0pt} \setlength{\itemsep}{0pt}}}
			\newcommand{\el}{\end{list}}
		
		\newtheorem{Thm}{Theorem}[section]
		\newtheorem{Lem}[Thm]{Lemma}
		
		\newtheorem{Def}[Thm]{Definition}
		\newtheorem{Exm}[Thm]{Example}
		\newtheorem{Rem}[Thm]{Remark}
		\newtheorem{Cor}[Thm]{Corollary}
		
        \newtheorem{Nota}[Thm]{Notation}

		\title{}
		\author{}
		\date{}
		\usepackage{amssymb}

		\usepackage{hyperref}
		\hypersetup{
			colorlinks,
			citecolor=blue,
			filecolor=black,
			linkcolor=blue,
			urlcolor=black
		}
\begin{document}
	
	\title{Cohomology of Lie algebroids over Topological ringed spaces}
	\author{ Abhishek Sarkar}

	\begin{abstract}
		We consider Lie algebroids over a topological ringed space as quasicoherent sheaves of Lie-Rinehart algebras.
		We express hypercohomology  for a locally free Lie algebroid (not necessarily of finite rank) as a derived functor, and simplify it via \v{C}ech cohomology.
		Furthermore, we define the Hochschild hypercohomology of a  sheaf of generalized bialgebras (using a derived functor) and  study the cases of the universal enveloping algebroid and of the jet algebroid  of a Lie algebroid. 
		In the sequel, we present a version of Hochschild-Kostant-Rosenberg theorem for a locally free Lie algebroid, as well as its dual version.
	\end{abstract}
	\footnote{AMS Mathematics Subject Classification : $16$S$30$, $16$T$10, $$17$B$56$, $32$L$10$, $53$D$17$.
		
	}
	\keywords{ Lie algebroids, Lie-Rinehart algebras, logarithmic derivations, (logarithmic) tangent sheaves,  universal enveloping algebras, bialgebras, Lie algebroid cohomology, Hochschild cohomology}
	\maketitle
	\section{Introduction}
	The notion of \emph{Lie algebroids} plays a prominent role in differential geometry and mathematical physics as they represent generalized infinitesimal symmetries of spaces, which are related to the corresponding global symmetries of spaces described by Lie groupoids \cite{KM}. Lie algebroids over a $C^{\infty}$-manifold $M$ are joint generalization of the tangent vector bundle $TM$ of $M$ and Lie algebras \cite{MM}. 
An algebraic analogue of Lie algebroids, known as \emph{Lie-Rinehart algebras}, is used to study more general situations \cite{JH,MK,LM,UB}.
    	In the context of complex geometry and algebraic geometry, B. Pym, A. Sarkar et.al. studied Lie algebroids over analytic spaces \cite{BP, AS} and M. Kapranov, D. Calaque, M. Van den Bergh, U. Bruzzo et.al. studied Lie algebroids over algebraic varieties (or schemes) \cite{MK,CV,UB} in the \emph{sheaf theoretic language}. These Lie algebroids (over such a space) are joint generalization of the \emph{tangent sheaf} and \emph{sheaves of Lie algebras} \cite{MK,SR,BP}. 
    \vspace{.05 cm}

    We study \emph{Lie algebroids over a topological ringed space}  $(X, \mathcal{O}_X)$ as certain \emph{quasicoherent sheaves of Lie-Rinehart algebras} \cite{MK,JV, AA, PS}, where $\mathcal{O}_X $ is a sub-sheaf of algebras of the sheaf of continuous functions $C_X^0$ on $X$ (in short, we refer $(X, \mathcal{O}_X)$ as an $a$-space) \cite{SR,AM}. 
	 This framework unifies the concept of Lie algebroids across the three types of base spaces as special cases. To study calculus on the smooth and singular geometric objects in a consistent manner, we need the notion of Lie algebroids within the algebro-geometric settings.
  It allows one to treat several geometric structures, such as Poisson analytic spaces \cite{BP}, singular foliations or  involutive distributions  \cite{RF,BP,LG}, as well as free Lie algebroids \cite{MK}. 
	In particular, the \emph{sheaf of logarithmic derivations} for some (principal or free) divisor of a complex manifold (or a smooth algebraic variety)  \cite{CMD,BP,AA}, and the path algebroid of a smooth manifold \cite{MK}, are key object of study in this context.
	
	 \vspace{.05 cm}

	Associated with a Lie algebroid over an $a$-space, there are two canonical sheaves of generalized bialgebras: the \emph{universal enveloping algebroid} and the \emph{jet algebroid} \cite{MK,CV,CRV,AP,BP,AA,SA}. Both play crucial roles in studying the homological algebra of a Lie algebroid. Here, the \emph{sheaves of generalized bialgebras} serve as a sheaf-theoretic analogue (or global version) of the notion of 
	$R/{\mathbb{K}}$-bialgebras \cite{MM} or left $R$-bialgebroids \cite{KP}, referred to as $\mathcal{O}_X/{\mathbb{K}_X}$-bialgebras \cite{AA}.
	The universal enveloping algebroid $\mathscr{U}(\mathcal{O}_X,\mathcal{L})$  of a Lie algebroid $\mathcal{L}$ generalizes the notion of \emph{sheaf of differential operators} on a manifold \cite{TS,AA,SR,SA}. It is \emph{sheafification of the presheaf} of universal enveloping algebras of Lie-Rinehart algebras associated with each space of sections of $\mathcal{L}$ \cite{UB,AA}. It has a canonical $\mathcal{O}_X/{\mathbb{K}_X}$-bialgebra structure \cite{AA}, similar structures are present in \cite{MK,CV} by different names.
	Moreover the dual of $\mathscr{U}(\mathcal{O}_X,\mathcal{L})$, i.e. $\mathscr{J}(\mathcal{O}_X,\mathcal{L}) := \mathscr{H}om_{\mathcal{O}_X}(\mathscr{U}(\mathcal{O}_X,\mathcal{L}),\mathcal{O}_X)$ is the jet algebroid of $\mathcal{L}$, generalizes the notion of sheaf of jets on a manifold (see \cite{MK,CV,DRV,BP,AP}). 
	It is sheafification of the presheaf of jet algebras associated with the sheaf of Lie-Rinehart algebras \cite{AP,AA}. 
It also has a canonical $\mathcal{O}_X/{\mathbb{K}_X}$-bialgebra structure induced from the structure of $\mathscr{U}(\mathcal{O}_X,\mathcal{L})$ \cite{AA}. On the other hand, an $\mathcal{O}_X/{\mathbb{K}_X}$-bialgebra gives rise to a Lie algebroid over the $a$-space $(X, \mathcal{O}_X)$, consisting of its sheaf of primitive elements \cite{AA}.
	 \vspace{.05 cm}

	Cohomology of a Lie algebroid generalizes both Lie algebra cohomology (Chevalley-Eilenberg cohomology) and de Rham cohomology, and is commonly refereed to as \emph{Chevalley-Eilenberg-de Rham cohomology} \cite{KM,BP,AP,BRT,UB}. 
The cohomology of a Lie algebroid $\mathcal{L}$, where $\mathcal{L}$ is a locally free $\mathcal{O}_X$-module of finite rank, over a Noetherian separated scheme  or over a complex manifold $(X,\mathcal{O}_X)$, is described using the derived functor $Ext$ by Bruzzo \cite{UB} and Pym \cite{BP}.
	In non-commutative geometry \cite{II,VG,TS}, the dual pair consisting of the Hochschild cohomology and homology $(HH^\bullet(A), HH_{\bullet}(A))$ of an associative $\mathbb{K}$-algebra $A$ serves as a non-commutative analogue to the classical dual pair $(\wedge^\bullet_{\mathcal{O}_X} \mathcal{T}_X, \wedge^\bullet_{\mathcal{O}_X} \Omega^1_X)$, which consists of the \emph{sheaf of multivector fields} and the \emph{sheaf of differential forms} over a smooth $a$-space $(X, \mathcal{O}_X)$. This correspondence is given by the standard \emph{Hochschild-Kostant-Rosenberg} (shortly, HKR) theorem. In \cite{CV, DRV, CRV}, a version of the HKR theorem is established by D. Calaque, C. A. Rossi, M. Van den Bergh, using the dual pair  $\mathscr{U}(\mathcal{O}_X,\mathcal{L})$ and $\mathscr{J}(\mathcal{O}_X,\mathcal{L})$, to explore precalculus up to homotopy for locally free Lie algebroids of finite rank. For that, an isomorphism of Gerstenhaber algebras is introduced, given by a twisted version of the HKR morphism. In this article, we deduce analogous cohomological results using a different approach to include more general situations, such as relaxing the \emph{finite rank condition}. We don't require the twisted version of the HKR morphism, as we work on the level of graded vector space isomorphism, rather that of Gerstenhaber algebra isomorphism.
    To derive such results, we present a \emph{global version} of cohomological results on Lie-Rinehart algebras \cite{GR,JH,KP,HOM}.

 \vspace{.05 cm}
	In Section \ref{Sec 2}, we recall some of the preliminary notions associated with a Lie algebroid that are required for the article.
	In Section \ref{Sec 3}, we consider \emph{Lie algebroid hypercohomology} over certain special $a$-spaces.
	Here, we show that one can define Lie algebroid cohomology, using a slightly different cochain complex of sheaves of quasicoherent  $\mathcal{O}_X$-modules.
	 The cohomology of a locally free Lie algebroid $\mathcal{L}$  over an $a$-space $(X, \mathcal{O}_X)$ with coefficient in an $\mathcal{L}$-module $\mathcal{E}$, is given by	$\mathbb{H}^\bullet(\mathcal{L}, \mathcal{E}) \cong Ext^\bullet_{\mathscr{U}(\mathcal{O}_X,\mathcal{L})}(\mathcal{O}_X,\mathcal{E})$ (where $\mathcal{L}$ is not necessarily of finite rank $\mathcal{O}_X$-module).
Here, we use the result of Lie-Rinehart cohomology as derived functor for a projective Lie-Rinehart algebra given by G. Rinehart \cite{GR,JH} as a local description (see Theorem \ref{Lie algebroid Cohomology as derived functor}). We apply this result for the sheaf of logarithmic derivations and for $\mathcal{O}_X/{\mathbb{K}_X}$-bialgebras. Furthermore, we express the Lie algebroid cohomology in terms of \emph{\v{C}ech cohomology} (see Section \ref{Cech cohomology derived}) by considering a good open cover \cite{SR,MS,BRT}. As an application of the cohomology groups, we consider \emph{abelian Lie algebroid extensions} (see Theorem \ref{extensions of L by E}) and (first) \emph{Chern classes} (see Section \ref{first Chern class}) in this general set up, extending results of H. Makestaad \cite{HM,HOM}.
	In Section \ref{Sec 4}, we consider algebraic (analytic) de Rham cohomologies for some free divisors associated with principal ideal sheaves \cite{AA} and compute the corresponding \emph{logarithmic de Rham cohomologies} \cite{CMD,Mor}. We view these cohomologies as Lie algebroid (hyper)cohomology  \cite{BP,BRT,UB}.

     \vspace{.05 cm}

	In the first part of Section \ref{Sec 5}, we define \emph{Hochschild hypercohomology} of an $\mathcal{O}_X/{\mathbb{K}_X}$-bialgebra. In particular, we study Hochschild hypercohomology of $\mathscr{U}(\mathcal{O}_X,\mathcal{L})$ and of $\mathscr{J}(\mathcal{O}_X,\mathcal{L})$ associated with a Lie algebroid $\mathcal{L}$.
	After that, we present a version of Hochschild-Kostant-Rosenbergh (HKR) theorem
	for locally free Lie algebroids (locally free $\mathcal{O}_X$-modules but not necessarily of finite rank) over any of the special $a$-spaces (see Section \ref{Hochschild coho of universal enveloping alg}). It provides an isomorphism of graded vector spaces between Hochschild hypercohomology of $\mathscr{U}(\mathcal{O}_X,\mathcal{L})$ and hypercohomology of the \emph{sheaf of $\mathcal{L}$-poly vector fields}, that is $\mathbb{H}H^\bullet(\mathscr{U}(\mathcal{O}_X,\mathcal{L})) \cong \mathbb{H}^\bullet(X, \wedge^\bullet_{\mathcal{O}_X} \mathcal{L})$ (see Theorem \ref{Generalized HKR thm}). We obtain this result by locally using the algebraic counterpart for projective Lie-Rinehart algebras established  by N. Kowalzig and H. Posthuma in \cite{KP}. 
	We derive a result in sheaf cohomology context (see Lemma \ref{hyper & sheaf}) to simplify the hypercohomology $\mathbb{H}^\bullet(X, \wedge^\bullet_{\mathcal{O}_X} \mathcal{L})$ (see Corollary \ref{simplify HKR}). Moreover, we present the HKR theorem in some of the special cases. Next we discuss about the dual version of HKR theorem in the generalize setup (see Section \ref{Hochschild coho of jet algebroid}). Thus we show that for a locally free Lie algebroids $\mathcal{L}$ of finite rank over any of the special $a$-space $(X, \mathcal{O}_X)$, we get a canonical isomorphism of vector spaces between Hochschild hypercohomology of $\mathscr{J}(\mathcal{O}_X, \mathcal{L})$ and the Lie algebroid (hyper)cohomology of $\mathcal{L}$, that is $\mathbb{H}H^\bullet(\mathscr{J}(\mathcal{O}_X,\mathcal{L})) \cong \mathbb{H}^\bullet(\mathcal{L}, \mathcal{O}_X)$ (See Theorem \ref{Generalized dual HKR thm}). This result we get by locally using the algebraic counterpart  as described in \cite{KP}.
	Both of these hypercohomologies are obtained by the derived functor $Cotor$. At last, we apply the dual HKR theorem for the tangent sheaf $\mathcal{L}:=\mathcal{T}_X$ over non-singular $a$-spaces and obtain some  results using the (smooth, analytic or algebraic) de Rham cohomology \cite{VID,MS} of $X$.

\vspace{.05 cm}

This article is an updated version of \emph{the authors' unpublished work} \cite{AS1}, with certain modifications.

	\section{Lie algebroids over an $a$-space $(X, \mathcal{O}_X)$}  \label{Sec 2}
	In this section, we recall some standard notions, namely  topological ringed spaces, Lie algebroids, universal enveloping algebroids,  generalized bialgebras etc in algebro-geometric language (see \cite{AA,PS}). Then, we recall some important relationships among them. In the later sections, we use these notions to consider cohomology theoretic results for such Lie algebroids.

	The notation $\mathbb{K}$ is used for $\mathbb{R, C}$ (the real or complex number fields respectively) or an algebraically closed field of characteristic zero. We denote the constant sheaf by $\mathbb{K}_X$ on a topological space $X$ with stalk $\mathbb{K}$. The sheaf of $\mathbb{K}$-valued (with one of the standard topologies) continuous functions on $X$ is denoted by $C_X^0$.
	\subsection{Lie algebroids over topological ringed spaces} \label{Lie algebroids}
	Here we consider some special locally ringed spaces \cite{CW}, and then define Lie algebroids over these spaces \cite{PS}.
	\begin{Def}
		Let $X$ be a topological space and $\mathcal{O}_X$ a $\mathbb{K}_X$-subalgebra of the sheaf of continuous functions $C_X^0$ on $X$. 
		The pair $(X, \mathcal{O}_X)$ is said to be a topological ringed space, or simply, an  $a$-space. 
	\end{Def}	
     For an $a$-space $(X, \mathcal{O}_X)$, the sheaf of $\mathbb{K}_X$-linear derivations of $\mathcal{O}_X$ is denoted by $\mathcal{D}er_{\mathbb{K}_X}(\mathcal{O}_X)$.
     This is the subsheaf of the sheaf endomorphisms $\mathcal{E}nd_{\mathbb{K}_X}(\mathcal{O}_X)$ on $\mathcal{O}_X$ satisfies the Leibniz rule. More precisely, for an open set $U$ of $X$, $\mathcal{D}er_{\mathbb{K}_X}(\mathcal{O}_X)(U):= \{D: \mathcal{O}_U \rightarrow \mathcal{O}_U~~\mbox{ is a derivation}\} \subset \mathcal{E}nd_{\mathbb{K}_X}(\mathcal{O}_U)$, where $\mathcal{O}_U=\mathcal{O}_X|_U$. That is, 
        a section $D$ of $\mathcal{D}er_{\mathbb{K}_X}(\mathcal{O}_X)$ on $U$ is a $\mathbb{K}$-linear derivation $D$ of $\mathcal{O}_X(U)$ (or, $D \in Der_{\mathbb{K}}(\mathcal{O}_X(U)))$ such that
        for open subsets $W, V$ of $U$ with $W \subset V$, satisfies the compatibility condition: 
        $$res_{VW}^{\mathcal{O}_X} \circ D_V= D_W \circ res_{VW}^{\mathcal{O}_X},$$
where $D_{V'}= D|_{V'}$, i.e.  $D_{V'}=res_{UV'}^{\mathcal{E}nd_{\mathbb{K}_X}(\mathcal{O}_X)} (D)$, for an open subset $V'$ of $U$.

         The sheaf $\mathcal{D}er_{\mathbb{K}_X}(\mathcal{O}_X)$ over such a space $(X, \mathcal{O}_X)$ has a canonical $ \mathcal{O}_X$-module and $\mathbb{K}_X$-Lie algebra structure, which are compatible via a Leibniz rule \cite{SR,MK,BP}.  
  
	\begin{Rem}\label {special $a$-spaces} 
		Consider smooth manifolds, complex manifolds, analytic spaces and algebraic varieties, all with their associated structure sheaf, as special $a$-spaces. For a special $a$-space $(X, \mathcal{O}_X)$, the sheaf $\mathcal{D}er_{\mathbb{K}_X}(\mathcal{O}_X)$ is known as the tangent sheaf of $X$, and is denoted by $\mathcal{T}_X$. 
        The tangent sheaf of $(X, \mathcal{O}_X)$ has the standard compatible $\mathbb{K}_X$-Lie algebra and a $($quasi$)$coherent $\mathcal{O}_X$-module structure \cite{MK,BP}. 
        
        A special $a$-space $(X, \mathcal{O}_X)$ is  non-singular (smooth) if and only if its tangent sheaf is a locally free $\mathcal{O}_X$-module of finite rank.
	\end{Rem}
     Recall that, a $(\mathbb{K}, A)$-\emph{Lie-Rinehart algebra} $L$ is a $\mathbb{K}$-Lie algebra  $L$ with an $A$-module structure (where $A$ is a commutative unital $\mathbb{K}$-algebra) and an $A$-module homomorphism $\rho:L\rightarrow Der_{\mathbb{K}}(A)$ such that 
\begin{itemize}
\item the map $\rho$ is also a $\mathbb{K}$-Lie algebra homomorphism,  
\item $ [x, ry]=r[x,y]+\rho(x)(r)y~~\mbox{ for all}~ x,y\in L,~r\in A.$
\end{itemize}

 The standard notion of Lie algebroids is defined on a smooth (real) vector bundle over a smooth manifold. The space of global sections of a Lie algebroid forms a Lie-Rinehart algebra, capturing the entire structure.
 In contrast, for holomorphic (or algebraic) vector bundles over complex manifolds (or algebraic varieties), the space of global sections does not capture the whole information. A sheaf-theoretic approach is necessary for analytic and algebraic geometry, applicable to both singular and non-singular cases.
	
	\begin{Def}  A Lie algebroid $\mathcal{L}$ over an $a$-space $(X,\mathcal{O}_X)$ is a quasicoherent sheaf of $(\mathbb{K}_X, \mathcal{O}_X)$-Lie-Rinehart algebras.
		That is, $\mathcal{L}$ is  a quasicoherent $\mathcal{O}_X$-module with a $\mathbb{K}_X$-Lie algebra structure on it defined by a bracket $[\cdot,\cdot]$, and
        equipped with a homomorphism  $\mathfrak{a}:(\mathcal{L},[\cdot,\cdot])\rightarrow (\mathcal{D}er_{\mathbb{K}_X}(\mathcal{O}_X),[\cdot,\cdot]_c)$ of $\mathcal{O}_X$-modules and $\mathbb{K}_X$-Lie algebras, called the anchor map. The map $\mathfrak{a}$  satisfies the Leibniz rule:
		$$[D,f D'] = f[D,D']+\mathfrak{a}(D)(f) D',$$ for $f \in \mathcal{O}_X$ and $D,D' \in \mathcal{L}$.
		We denote this Lie algebroid as $(\mathcal{L}, [\cdot, \cdot], \mathfrak{a})$ or simply by $\mathcal{L}$.
	\end{Def}
	\begin{Rem}
		Lie algebroids over a smooth manifold $X$ is equivalent to locally free Lie algebroids of finite rank over the $a$-space $(X, C_X^\infty)$, where $ C_X^\infty$ is the sheaf of $C^\infty$-functions on $X$.
	\end{Rem}
	\emph{Lie algebroid morphisms.} A homomorphism of Lie algebroids over an $a$-space $(X, \mathcal{O}_X)$ 
	$$\phi :  (\mathcal{L}_1,[\cdot,\cdot]_1,\mathfrak{a}_1)\rightarrow (\mathcal{L}_2,[\cdot,\cdot]_2,\mathfrak{a}_2)$$ is a sheaf homomorphism of Lie-Rinehart algebras, i.e.
	\begin{itemize}
		\item $\phi: \mathcal{L}_1 \rightarrow \mathcal{L}_2$ is an $\mathcal{O}_X$-module homomorphism,
		\item $\phi :  (\mathcal{L}_1,[\cdot,\cdot]_1)\rightarrow (\mathcal{L}_2,[\cdot,\cdot]_2)$ is a $\mathbb{K}_X$-Lie algebra homomorphism,
		\item compatibility condition: $\mathfrak{a}_2 \circ \phi= \mathfrak{a}_1$.
	\end{itemize}
    \emph{Locally free Lie algebroids.}  A Lie algebroid over $(X, \mathcal{O}_X)$ is called a locally free Lie algebroid (of finite rank) if its underlying $\mathcal{O}_X$-module is locally free (of finite rank). Similarly, we refer to a Lie algebroid as (quasi)coherent Lie algebroid if its underlying $\mathcal{O}_X$-module structure is (quasi)coherent (see \cite{MK}).

For a (smooth or holomorphic) vector bundle $E$ over a (smooth or complex) manifold $X$, the sheaf of sections of $E$ over $X$ is denoted by $\Gamma_X(E)$. Here, $\Gamma_X$ represents the sheaf of sections functor.

	\begin{Exm} The standard Lie algebroid structure on the tangent sheaf of an $a$-space $(X, \mathcal{O}_X)$ is described by the pair $(\mathcal{D}er_{\mathbb{K}_X}(\mathcal{O}_X), [\cdot,\cdot]_c, id)=:\mathcal{T}_X$. In particular, when $X$ is a real smooth manifold $($complex manifold$)$ with the structure sheaf $\mathcal{O}_X$, the sheaf of smooth $($holomorphic$)$ vector fields $\Gamma_X(TX)$ on $X$ is isomorphic to  $\mathcal{D}er_{\mathbb{K}_X}(\mathcal{O}_X)$, forming a locally free Lie algebroid of rank equals to dimension of the manifold  $X$.
	\end{Exm}	
    \begin{Exm} \label{Cotantgent Poisson}
		 The cotangent sheaf $($i.e. the sheaf of differential $1$-forms$)$ $\Omega^1_X$ over a Poisson manifold $X$ with Poisson bi-vector field $\pi \in H^0(X, \wedge^2 \mathcal{T}_X)$, has a canonical Lie algebroid structure, described as follows.
        The anchor map $\tilde{\pi}: \Omega^1_X \rightarrow \mathcal{T}_X$ is defined by $\tilde{\pi}(\omega) (f)= \pi(\omega, df)$ for $f \in \mathcal{O}_X$, and the Lie bracket on $\Omega^1_X$ is defined as 
        $[\omega, \omega']:= L_{\tilde{\pi}(\omega)} \omega' - L_{\tilde{\pi}(\omega')} \omega - d(\pi(\omega, \omega'))$, for $f \in \mathcal{O}_X$ and $\omega, \omega' \in \Omega^1_X$ (see \cite{RF, BP}).
	\end{Exm}	
	\begin{Exm} \label{Foliation}
		A singular foliation  $\mathcal{F}$ on a  real smooth manifold or a complex manifold $(X, \mathcal{O}_X)$ is an $\mathcal{O}_X$-submodule of the Lie algebroid $\mathcal{T}_X$ $($or $\mathfrak{X}_X)$, which is $(i)$ stable under the Lie bracket and $(ii)$ locally finitely generated. It provides a Stefan-Sussmann or Nagano distribution on $X$, forms a Lie algebroid over $(X, \mathcal{O}_X)$.  
		In particular, regular foliations arise from involutive subbundles of the tangent bundle, providing Frobenius distributions $($see \cite{RF, BP, LG} for the smooth and holomorphic cases$)$. 
	\end{Exm}
	\begin{Exm}
		\label{analytic spaces} 
		Let $X$ be a complex manifold and $\mathcal{O}_X$ be the sheaf of holomorphic functions. The vanishing set $($or the zero locus$)$ $Y:= V(\mathcal{I})$ of a ideal-sheaf $\mathcal{I} \subset \mathcal{O}_X$ is a subspace of $X$ which is not necessarily a submanifold $($it may have singularity as well$)$. The sheaf of functions $\mathcal{O}_Y := \mathcal{O}_X/{\mathcal{I}}$ on $Y$ is its structure sheaf and the pair $(Y, \mathcal{O}_Y)$ is an analytic space (see in \cite{BP}).

		Here, the tangent sheaf $\mathcal{T}_Y := (\mathcal{D}er_{\mathbb{C}_Y}(\mathcal{O}_Y),  [\cdot,\cdot]_c, id_Y)$ is not necessarily a locally free Lie algebroid over $(Y, \mathcal{O}_Y)$.
		Consider the sheaf of logarithmic derivations $($forms a generalized involutive distribution or singular foliation$)$ \cite{BP, AA} as 
		$$\mathcal{T}_X(-log~ Y) := \{D \in \mathcal{T}_X : D(\mathcal{I}) \subset \mathcal{I}\} \hookrightarrow \mathcal{T}_X$$
		$($geometrically, it represents the sheaf of vector fields on $X$ that are tangent to $Y$ for a smooth divisor $Y)$ with the canonical Lie algebroid structure. It is associated with $\mathcal{T}_Y$ via the canonical map 
		$$\rho : \mathcal{T}_X(-log~ Y) \rightarrow \mathcal{T}_Y$$
        defined by $\rho (g D) = \bar{g} \tilde{D}$, for $g \in \mathcal{O}_X$, $D \in\mathcal{T}_X(-log~Y)$ as $\tilde{D}(\bar{f}) = \bar{D(f)}$, where $\bar{h}=h + \mathcal{I} \in \mathcal{O}_Y$ for $h \in \mathcal{O}_X$.
	\end{Exm}
	
	\subsection{Universal enveloping algebroid of a Lie algebroid}
	Let $(\mathcal{L},[\cdot,\cdot],\mathfrak{a})$ be  a Lie algebroid  over any of the special $a$-space $(X, \mathcal{O}_X)$.  For each open set $U$ of $X$, we find  the universal enveloping algebra $\mathcal{U}(\mathcal{O}_X(U),\mathcal{L}(U))$ of the $(\mathbb{K},\mathcal{O}_X(U))$-Lie-Rinehart algebra $\mathcal{L}(U)$ or of the Lie-Rinehart pair $(\mathcal{O}_X(U),\mathcal{L}(U))$ (see \cite{GR,JH,MM}). The sheafification of the cannonical presheaf: 
	$$U \mapsto \mathcal{U}(\mathcal{O}_X(U),\mathcal{L}(U)),$$ is known as the universal enveloping algebroid of the Lie algebroid $\mathcal{L}$, and denoted by $\mathscr{U}(\mathcal{O}_X, \mathcal{L})$ (see \cite{DRV,BP,UB,TS,AA} for details).
	
	From the construction of $\mathscr{U}(\mathcal{O}_X, \mathcal{L})$, it is an associative $\mathbb{K}_X$-algebra and $\mathcal{O}_X$-bimodule. Moreover,
	there is a canonical  $\mathbb{K}_X$-algebra monomorphism $\iota_X: \mathcal{O}_X \hookrightarrow \mathscr{U}(\mathcal{O}_X,\mathcal{L})$ and an $\mathcal{O}_X$-linear map $\iota_{\mathcal{L}}: \mathcal{L} \rightarrow \mathscr{U}(\mathcal{O}_X,\mathcal{L})$. 
	The associative $\mathbb{K}_X$-algebra $\mathscr{U}(\mathcal{O}_X, \mathcal{L})$ is generated by $\mathcal{O}_X$ and $\iota_{\mathcal{L}}(\mathcal{L})$ satisfy the following identities:
	\begin{center}
		$\bar{D} ~ \bar{D'} - \bar{D'} ~ \bar{D}= \overline{[D,~D']}$,~ $ \bar{D} ~ f - f ~ \bar{D}= \mathfrak{a}(D)(f)$,
	\end{center}
	where $D, D' \in \mathcal{L}$, ~$f \in \mathcal{O}_X$, and $\bar{D}= \iota_{\mathcal{L}}(D)$ for all $D\in \mathcal{L}$.
	
	Hence, the map $\iota_{\mathcal{L}}$ can also be viewed as a $\mathbb{K}_X$-Lie algebra homomorphism.
	\begin{Rem}\label{Uni prop}
		The universal enveloping algebroid $\mathscr{U}(\mathcal{O}_X, \mathcal{L})$ of $\mathcal{L}$ is characterized by the following universal property: \label{Universal}
		
		Let $\mathcal{A}$ be a sheaf of unital associative $\mathbb{K}_X$-algebra with sheaf homomorphisms $\phi: \mathcal{O}_X \rightarrow \mathcal{A}$ of $\mathbb{K}_X$-unital algebras and $\psi: (\mathcal{L}, [\cdot,\cdot]) \rightarrow (\mathcal{A}, [\cdot,\cdot]_c)$ of  $\mathcal{O}_X$-linear $\mathbb{K}_X$-Lie algebras such that $\phi(f)\psi(D) = \psi(fD)$ and $[\psi(D), \phi(f)]_c = \phi(\mathfrak{a}(D(f)))$ holds for $f \in \mathcal{O}_X$ and $D \in \mathcal{L}$. Then, there exists a unique $\mathcal{O}_X$-linear homomorphism of unital $\mathbb{K}_X$-algebras $\widetilde{\psi} : \mathscr{U}(\mathcal{O}_X, \mathcal{L}) \rightarrow \mathcal{A}$ such that $\widetilde{\psi} \circ \iota_X = \phi$ and $\widetilde{\psi} \circ \iota_{\mathcal{L}} = \psi$.
	\end{Rem}
	
	\begin{Rem}\label{embed}
		For a locally free Lie algebroid $\mathcal{L}$ over an $a$-space $(X, \mathcal{O}_X)$, the homomorphism $\iota_{\mathcal{L}}: \mathcal{L} \rightarrow \mathscr{U}(\mathcal{O}_X,\mathcal{L})$ become an embedding (as $\mathbb{K}_X$-algebras and $\mathcal{O}_X$-modules). 
	\end{Rem}
	\begin{Exm}
		For a non-singular special $a$-space $(X,\mathcal{O}_X)$, the universal enveloping algebroid is isomorphic to the sheaf of differential operators $\mathcal{D}_X$ on $X$ (i.e. the sheaf of differential operators over $\mathcal{O}_X$, sometimes denoted as $\mathcal{D}iff(\mathcal{O}_X))$ \cite{MK, TS}, i.e. $$\mathscr{U}(\mathcal{O}_X, \mathcal{T}_X) \cong \mathcal{D}iff(\mathcal{O}_X)=:\mathcal{D}_X.$$
	\end{Exm}
	\begin{Exm} \label{sheaf of log diff operators}
		The sheaf of logarithmic derivations  and sheaf of logarithmic differential operators are denoted as $\mathcal{T}_X(-log~ Y)$ and $\mathcal{D}_X(-log~ Y)$ respectively for a principal divisor $Y$ in some complex manifold or smooth algebraic variety $X$ \cite{AA}.
		
		In the case of a free divisor $Y$ in $X$ (i.e. $\mathcal{T}_X(-log~ Y)$ is a locally free $\mathcal{O}_X$-module \cite{CMD,BP,AA})  we have (sheafifying the local description given for the module of logarithmic derivations in \cite{LM}) the isomorphism
		$$\mathscr{U}(\mathcal{O}_X, \mathcal{T}_X(-log~ Y)) \cong  \mathcal{D}_X(-log~ Y).$$
	\end{Exm}
	\begin{Exm}\label{free Lie algebroid} In \cite{MK}, the notion of path algebroid $\mathcal{P}_X$ of a smooth manifold or smooth algebraic variety $X$ is constructed as the free Lie algebroid generated by tangent sheaf $\mathcal{T}_X$.  It forms a locally free sheaf of Lie-Rinehart algebras over $(X,\mathcal{O}_X)$ of infinite rank. The universal enveloping algebroid $\mathscr{U}(\mathcal{O}_X,\mathcal{P}_X) =:\mathbb{D}_X$ is described as sheaf of non commutative differential operators on $X$.
	\end{Exm}  
	\subsection{$\mathcal{O}_X/{\mathbb{K}_X}$-bialgebras} \label{general bialgebra}
	The notion $\mathcal{O}_X/{\mathbb{K}_X}$-bialgebras  \cite{AA} is considered as a sheaf theoretic analogue of the notion $R/{\mathbb{K}}$-bialgebras \cite{MM} or left bialgebroids \cite{KP}. 
	As an example, we consider the universal enveloping algebroid $\mathscr{U}(\mathcal{O}_X,\mathcal{L})$ (and its dual) of a Lie algebroid $\mathcal{L}$ over $(X, \mathcal{O}_X)$.
	Similar structure is appeared in \cite{MK, CV, DRV, AP}, in complex and algebraic geometry context.
	\begin{Def} 
		A tuple $(\mathcal{A}, \Delta,\epsilon)$  is called $\mathcal{O}_X/{\mathbb{K}_X}$-bialgebra if the following conditions hold.

		\begin{itemize}
			\item $\mathcal{A}$ is a sheaf of unital associative $\mathbb{K}_X$-algebra extending the sheaf $\mathcal{O}_X$;
			\item  $\mathcal{A}$ is equipped with the morphism of $\mathcal{O}_X$-modules, the comultiplication 
			$\Delta : \mathcal{A}\rightarrow \mathcal{A}\otimes_{\mathcal{O}_X} \mathcal{A}$ and the counit $\epsilon : \mathcal{A}\rightarrow \mathcal{O}_X$;
			
			\item The image $\Delta(\mathcal{A})$ lies in a certain $\mathbb{K}_X$-subalgebra of $\mathcal{A} {\otimes}_{\mathcal{O}_X} \mathcal{A}$;
			\item  $\Delta(1) = 1\otimes 1$ and $\epsilon(1) = 1$, where $1$ is the unit of $ \mathcal{O}_X$;
			\item $\Delta(ab) = \Delta(a) \Delta(b)$,
			$\epsilon(ab) = \epsilon(a\epsilon(b))$ for any two sections $a,b$ in
			$\mathcal{A}$. 
		\end{itemize}
	\end{Def}
	We recall the canonical $\mathcal{O}_X/{\mathbb{K}_X}$-bialgebra structures on $\mathscr{U}(\mathcal{O}_X,\mathcal{L})$ and its dual in the following.

	\subsubsection{\textbf{The universal enveloping algebroid $\mathscr{U}(\mathcal{O}_X,\mathcal{L})$ as an $\mathcal{O}_X/{\mathbb{K}_X}$-bialgebra}}
	 \label{universal alg} 
	
	Let $(\mathcal{L},[\cdot,\cdot],\mathfrak{a})$ be a Lie algebroid over one of the special $a$-space $(X,\mathcal{O}_X)$. The sheaf $ \mathscr{U}(\mathcal{O}_X,\mathcal{L})$ has a canonical associative unital $\mathbb{K}_X$-algebras and an $\mathcal{O}_X$-bimodule structure as discussed above.
	
	Observe that $\mathscr{U}(\mathcal{O}_X,\mathcal{L})$ is a cocommutative counital $\mathcal{O}_X$-coalgebra in the sense that there is a natural cocommutative coassociative comultiplication  
	$\Delta : \mathscr{U}(\mathcal{O}_X,\mathcal{L}) \rightarrow \mathscr{U}(\mathcal{O}_X,\mathcal{L})\otimes_{\mathcal{O}_X} \mathscr{U}(\mathcal{O}_X,\mathcal{L})$, and a counit $\epsilon : \mathscr{U}(\mathcal{O}_X,\mathcal{L}) \rightarrow \mathcal{O}_X$ respectively,
	which are {locally} given by the following formulas \cite{DRV} 
	\begin{center}
		$\Delta(f) = f \otimes 1= 1 \otimes f$,\\
		$\Delta(\bar{D}) = \bar{D} \otimes 1 + 1 \otimes \bar{D}$,\\
		$\Delta (D' D'') = \sum D'_{(1)} D''_{(1)} \otimes D'_{(2)} D''_{(2)}$,\\  
		$\epsilon(\tilde{D}) = \tilde{D}(1)$, 
	\end{center}
	for a section $f$ of $\mathcal{O}_X$, for a section $D$ of $\mathcal{L}$ with $\bar{D}=\iota_{\mathcal{L}}(D)$, and for sections $D', D'',\tilde{D}$ of $\mathscr{U}(\mathcal{O}_X,\mathcal{L})$.
	Consequently, $(\mathscr{U}(\mathcal{O}_X,\mathcal{L}),\Delta,\epsilon)$ is a sheaf of $\mathcal{O}_X/{\mathbb{K}_X}$-bialgebras.
	
	\subsubsection{\textbf{The Jet enveloping algebroid $\mathscr{J}(\mathcal{O}_X,\mathcal{L})$ as an $\mathcal{O}_X/{\mathbb{K}_X}$-bialgebra}}
	 \label{jet alg} 
	The notion of jet algebroid of a Lie algebroid $\mathcal{L}$ over $(X, \mathcal{O}_X)$, is  defined as the dual of the universal enveloping algebroid:
	\begin{center}
		$\mathscr{J}(\mathcal{O}_X,\mathcal{L}) := \mathscr{H}om_{\mathcal{O}_X} (\mathscr{U}(\mathcal{O}_X,\mathcal{L}), \mathcal{O}_X) $.
	\end{center}
	For each open set $U$ of $X$, we find the jet algebras $\mathcal{J}(\mathcal{O}_X(U), \mathcal{L}(U))$ of the $(\mathbb{K},\mathcal{O}_X(U))$-Lie-Rinehart algebra $\mathcal{L}(U)$ (see \cite{CV,KP,AP}).
	The jet algebroid 	$\mathscr{J}(\mathcal{O}_X,\mathcal{L})$ is the sheafification of the canonical presheaf:
	$$ U \mapsto \mathcal{J}(\mathcal{O}_X(U), \mathcal{L}(U)).$$

One can dualize the structures on $\mathscr{U}(\mathcal{O}_X,\mathcal{L})$. The product on $\mathscr{J}(\mathcal{O}_X,\mathcal{L})$ is induced from the coproduct of $\mathscr{U}(\mathcal{O}_X,\mathcal{L})$ on each space of sections, which is locally defined by 
$$(\phi_1 \phi_2) (D) := \phi_1 (D_{(1)}) \phi_2 (D_{(2)})$$ for sections $\phi_1, \phi_2 \in \mathscr{J}(\mathcal{O}_X,\mathcal{L})$ and a section $D \in \mathscr{U}(\mathcal{O}_X,\mathcal{L})$.
By cocommutativity of $\mathscr{U}(\mathcal{O}_X,\mathcal{L})$, this defines a commutative algebra structure on $\mathscr{J}(\mathcal{O}_X,\mathcal{L})$. The unit for this multiplication is locally given by the left counit $\epsilon :\mathscr{U}(\mathcal{O}_X,\mathcal{L}) \rightarrow \mathcal{O}_X$, since
$$(\epsilon \phi) (D) = \epsilon (D_{(1)}) \phi (D_{(2)}) = \phi(\epsilon (D_{(1)}) D_{(2)})= \phi(D),$$
for a section $\phi$ in $\mathscr{J}(\mathcal{O}_X,\mathcal{L})$ and a section $D$ in $\mathscr{U}(\mathcal{O}_X,\mathcal{L})$.
The left and right $\mathcal{O}_X$-module structure on $\mathscr{U}(\mathcal{O}_X,\mathcal{L})$ provides $\mathcal{O}_X- \mathcal{O}_X$-bimodule structure on $\mathscr{J}(\mathcal{O}_X,\mathcal{L})$. 

The product on $\mathscr{U}(\mathcal{O}_X,\mathcal{L})$ descends to a sheaf homomorphism $m : \mathscr{U}(\mathcal{O}_X,\mathcal{L}) \otimes_{\mathbb{K}_X} \mathscr{U}(\mathcal{O}_X,\mathcal{L}) \rightarrow \mathscr{U}(\mathcal{O}_X,\mathcal{L})$. We can therefore dualize the product to obtain a coproduct $\Delta^* : \mathscr{J}(\mathcal{O}_X,\mathcal{L}) \rightarrow \mathscr{J}(\mathcal{O}_X,\mathcal{L}) \otimes_{\mathcal{O}_X} \mathscr{J}(\mathcal{O}_X,\mathcal{L})$ locally defined as
$$\phi(D D') =: \Delta^*(\phi)(D \otimes D') = \phi_{(1)}(D \phi_{(2)}(D')) $$

Associativity of the multiplication implies that $\Delta^*$ is coassociative. The counit for this coproduct is given by $\epsilon^* : \phi \mapsto \phi(1_{\mathscr{U}(\mathcal{O}_X,\mathcal{L})})$. It follows that $(\mathscr{J}(\mathcal{O}_X,\mathcal{L}), \Delta^*, \epsilon^*)$  is an $\mathcal{O}_X/{\mathbb{K}_X}$-bialgebra.

\begin{Rem}
	The jet algebroid is a formal groupoid that serves as the formal exponentiation of the Lie algebroid \cite{DRV}. In \cite{MK}, this structure is studied for the path algebroid $\mathcal{P}_X$ of a smooth variety $X$, referring to it as the formal path groupoid.
\end{Rem}
\begin{Rem} \label{primitives}
	For an $\mathcal{O}_X/{\mathbb{K}_X}$-bialgebra $(\mathcal{A}, \Delta,\epsilon)$, the sheaf of premitive elements is determined by the following assignment (see \cite{AA})
	$$U \mapsto \{a \in \mathcal{A}(U)~|~ \Delta(U)(a)= a \otimes 1 + 1 \otimes a\},$$
	denoted by $\mathscr{P}(\mathcal{A})$, provides a sheaf of $(\mathbb{K}_X, \mathcal{O}_X)$-Lie-Rinehart algebras as follows:
	$$\mathfrak{a} : (\mathscr{P}(\mathcal{A}), [\cdot,\cdot]_c)
	\rightarrow (\mathcal{D}er_{\mathbb{K}_X}(\mathcal{O}_X), [\cdot,\cdot]_c)$$ is defined by $\mathfrak{a}(D)(f) = \epsilon (D ~ f)$ for sections $f\in \mathcal{O}_X$ and $D\in \mathscr{P}({\mathcal{A}})$.
\end{Rem}
\begin{Rem}\label{cocomplete locally graded free bialgebras}For a smooth $($holomorphic$)$ Lie algebroid \cite{BRT} i.e. for a locally free sheaf of Lie-Rinehart algebras of finite rank over a smooth (complex) manifold, the universal enveloping algebroid is a cocomplete locally graded free $\mathcal{O}_X/{\mathbb{K}_X}$-bialgebra $($of finite type$)$ \cite{AA}.

	On the other hand, the tangent sheaf over the affine scheme $Spec (\mathbb{A}^{\mathbb{N}})$ $($where $\mathbb{A}^{\mathbb{N}}:= \mathbb{K}[x_i]_{i \in \mathbb{N}})$ \cite{VD,AA} and the free Lie algebroid $\mathcal{P}_X$  over a smooth manifold $X$ \cite{MK} are locally free sheaves of Lie-Rinehart algebras $($of infinite rank$)$. Thus, their universal enveloping algebroids are cocomplete locally graded free $\mathcal{O}_X/{\mathbb{K}_X}$-bialgebras $($of infinite type$)$ \cite{AA}.
    \end{Rem}
\section{Lie algebroid Cohomology over $a$-spaces} \label{Sec 3}

Lie algebroids over an $a$-space and homomorphisms among them are discussed in Section \ref{Lie algebroids}. These form a category, where kernel and image of a morphism provides new objects in that category.

We recall the notion of representations or modules of a Lie algebroid in this context. We use the fact that
the category of representations of a Lie algebroid $\mathcal{L}$ over $(X, \mathcal{O}_X)$ has enough injectives, since the category of left modules over the sheaf of rings $\mathscr{U}(\mathcal{O}_X, \mathcal{L})$ has enough injectives (see \cite{UB, MK, RH}).

The Chevalley-Eilenberg-de Rham cohomology (or, Lie algebroid cohomology) of a coherent Lie algebroid over a complex manifold or over a Noetherian separated scheme is expressed as derived functor $Ext$ \cite{BP,UB}.  In \cite{BP}, the result is stated for locally free Lie algebroid of finite rank over a complex manifold and in \cite{UB}, the result is proved for locally free Lie algebroid of finite rank over a Noetherian separated scheme. Here, we construct a more general cochain complex of sheaves, which coincides with the cochain complexes of sheaves used in the classical contexts \cite{KM,BP,BRT,UB}. For that we do not require the coherent (or locally finitely presented) condition on the underlying $\mathcal{O}_X$-module structure of a Lie algebroid.

Notably, for analytic and algebraic geometric setups, the base space $X$ does not need to be non-singular. However, in that case, the tangent sheaf, which becomes a coherent Lie algebroid, does not retain the property of being a locally free Lie algebroid $($see \cite{BP,AA}$)$.
 Moreover, we present the result in a simplified form by using \v{C}ech cohomology. As an application of Lie algebroid cohomology we consider abelian Lie algebroid extensions and Chern class of a module (locally free $\mathcal{O}_X$-module of finite rank) over the Lie algebroid.

To define the Lie algebroid cohomology of a Lie algebroid $\mathcal{L}$ over $(X, \mathcal{O}_X)$ with coefficients in an $\mathcal{O}_X$-module $\mathcal{E}$, we proceed as follows.

\subsubsection{Atiyah algebroid } For an (quasicoherent) $\mathcal{O}_X$-module $\mathcal{E}$, we form a Lie algebroid consisting of the sheaf of differential operators on $\mathcal{E}$ of order $\leq 1$ having scalar symbols \cite{MK,UB}, i.e. 
$$\mathcal{A}t(\mathcal{E}):= \{D\in \mathscr{E}nd_{\mathbb{K}_X}(\mathcal{E})~|~D(fs)=fD(s)+\sigma_D(f)s~~\mbox{ for a unique}~~  \sigma_D\in \mathcal{T}_X,~~\mbox{where}~~ f \in \mathcal{O}_X,~s\in \mathcal{E}\},$$
(thus, here $\sigma_D(f)=[D, f]_c \in \mathcal{O}_X$ for $D \in \mathcal{A}t(\mathcal{E})$ and $f \in \mathcal{O}_X$ holds) with the anchor map defined by 
\begin{center}
	$\sigma: \mathcal{A}t(\mathcal{E}) \rightarrow \mathcal{T}_X$ where $D \mapsto \sigma_D$
\end{center}
and the Lie bracket is commutator bracket. This Lie algebroid structure is so-called Atiyah algebroid of the $\mathcal{O}_X$-module $\mathcal{E}$. It provides a short exact sequence (s.e.s.) of Lie algebroids over $(X, \mathcal{O}_X)$
(an abelian Lie algebroid extension) as 
\begin{align}\label{Lie algebroid extension}
	0 \rightarrow \mathscr{E}nd_{\mathcal{O}_X}(\mathcal{E}) \hookrightarrow \mathcal{A}t(\mathcal{E})  \overset{\sigma}{\rightarrow} \mathcal{T}_X \rightarrow 0.
\end{align}

In particular, when $\mathcal{E}=\mathcal{O}_X$ for a non singular $a$-space $X$, we have $\mathcal{A}t(\mathcal{O}_X)\cong \mathcal{O}_X \oplus \mathcal{T}_X$, equals to the sheaf differential operators on $X$ of order $\leq 1$  with scalar symbol. The universal enveloping algebroid of the Atiyah algebroid $\mathcal{A}t(\mathcal{O}_X)$ is the sheaf of twisted differential operators over $X$ \cite{AP}.

A connection $\nabla : \mathcal{T}_X \rightarrow \mathcal{A}t(\mathcal{E})$ that satisfies $\sigma \circ \nabla=Id_{\mathcal{T}_X}$  provides a splitting of the s.e.s. (\ref{Lie algebroid extension}) as $\mathcal{O}_X$-modules. If its curvature is zero then the s.e.s. splits as Lie algebroids. More generally, for a Lie algebroid $(\mathcal{L}, [\cdot, \cdot], \mathfrak{a})$ on an $a$-space $(X, \mathcal{O}_X )$, an $\mathcal{L}$-connection on $\mathcal{E}$ is defined by an $\mathcal{O}_X$-linear map (see \cite{MK,UB})
\begin{align}\label{L-connection}
	\nabla: \mathcal{L}\rightarrow \mathcal{A}t(\mathcal{E})
\end{align}
$$\hspace{.5 cm}	D\mapsto \nabla_D$$
satisfying the Leibniz rule $\nabla_D(f~s)= f~\nabla_D(s)+\mathfrak{a}(D)(f)~s$ for  sections $f \in \mathcal{O}_X$, $D \in \mathcal{L}$ and $s \in \mathcal{E}$ $($to ensure compatibility, we choose $\sigma \circ \nabla= \mathfrak{a}$, i.e. $\sigma_{\nabla_{D}}=\mathfrak{a}(D)$ for $D\in \mathcal{L})$. 

Equivalently, it is described by a $\mathbb{K}_X$-linear map, if in addition $\mathcal{L}$ is a locally free $\mathcal{O}_X$-module of finite rank, as follows (see \cite{HM3, BP})
\begin{align} \label{L-connection special}
	\nabla : \mathcal{E} \rightarrow \Omega_{\mathcal{L}}^1 \otimes_{\mathcal{O}_X}\mathcal{E}
\end{align} 
satisfying the  Leibniz rule $\nabla(f~s)=f~\nabla s + {\mathfrak{a}^*}(df)\otimes s$, where $\Omega_{\mathcal{L}}^1:=\mathscr{H}om_{\mathcal{O}_X}(\mathcal{L}, \mathcal{O}_X)$ for a Lie algebroid $\mathcal{L}$ and $\Omega_{X}^1:=\Omega_{\mathcal{T}_{X}}^1$ (the Zarisky differentials \cite{CK}) together with ${\mathfrak{a}^*}:\Omega_X^1 \rightarrow \Omega_{\mathcal{L}}^1$ is the dual of the anchor map $\mathfrak{a}$.

The $\mathcal{L}$-connection on $\mathcal{E}$ is said to be flat if the map (\ref{L-connection}) is a Lie algebroid homomorphism  (i.e. the $\mathcal{L}$-curvature is zero), i.e. the map (\ref{L-connection}) additionally satisfies
\begin{align}\label{L-module}
	\nabla_{[D,D']}=[\nabla_D, \nabla_{D'}]_c
\end{align} 

In that case, $(\mathcal{E},\nabla)$ is said to be a representation of $\mathcal{L}$ or is called an $\mathcal{L}$-module.

Next, we consider some special cases where $X$ is non-singular $a$-spaces:

$(i)$ For $\mathcal{L}=\mathcal{T}_X$ and $\mathcal{E}$ is a locally free $\mathcal{O}_X$-module of finite rank, there is a $\mathcal{T}_X$-connection on $\mathcal{E}$ always exists if $X$ either a smooth manifold, or a Stein manifold, or an affine variety (see \cite{MA}).

$(ii)$ For $\mathcal{L}=\mathcal{T}_X(- log~ Y)$, if a connection exists on an $\mathcal{O}_X$-module $\mathcal{E}$, it is called a logarithmic connection on $\mathcal{E}$. It forms a meromorphic connection with simple poles along the divisor $Y$ (see \cite{CMD,BP}).

If $(\mathcal{E}, \nabla)$ is an $\mathcal{L}$-module, then the Lie algebroid morphism
\begin{center}
	$\nabla : \mathcal{L} \rightarrow \mathcal{A}t(\mathcal{E}) $
\end{center}  
extends to an $\mathcal{O}_X$-linear  homomorphism of $\mathbb{K}_X$-algebras (using Remark \ref{Universal})
$$\tilde{\nabla} : \mathscr{U}(\mathcal{O}_X,\mathcal{L}) \rightarrow \mathscr{E}nd_{\mathbb{K}_X}(\mathcal{E}),$$
making $\mathcal{E}$ into a left $\mathscr{U}(\mathcal{O}_X,\mathcal{L})$-module. If $\mathcal{L}$ is locally free $\mathcal{O}_X$-module and $\mathcal{E}$ be a $\mathscr{U}(\mathcal{O}_X,\mathcal{L})$-module, then the restriction of the action of $\mathscr{U}(\mathcal{O}_X,\mathcal{L})$ to $\mathcal{L}$ (by using the canonical embedding of $\mathcal{L}$ in $\mathscr{U}(\mathcal{O}_X,\mathcal{L})$ as described in Remark \ref{embed}), 
provides an $\mathcal{L}$-module structure on $\mathcal{E}$. In this case, the category of $\mathcal{L}$-modules and the category of left $\mathscr{U}(\mathcal{O}_X,\mathcal{L})$-modules are equivalent. It helps in studying  homological algebra with $\mathcal{L}$-modules (here $\mathcal{L}$ is equipped with a non-associative $\mathbb{K}_X$-algebra structure) by treating them as modules over the associative $\mathbb{K}_X$-algebra $\mathscr{U}(\mathcal{O}_X,\mathcal{L})$.

\subsubsection{Hypercohomology}  \label{Hypercohomology} To recall the notion of hypercohomology (see \cite{VID, UB}), we first fix some notions.

Let $\mathcal{O}$ be a sheaf of associative $\mathbb{K}$-algebras over a topological space $X$ with $\mathcal{O}$-modules $\mathcal{C}^i$, for $i \geq 0$. Consider a cochain complex $(\mathcal{C}^\bullet, d)$ where $\mathcal{C}^\bullet:= \oplus_i \mathcal{C}^i$ is a graded $\mathcal{O}$-module and $d:\mathcal{C}^\bullet \rightarrow \mathcal{C}^{\bullet +1}$ is the co-boundary map, i.e. $d$ is a $\mathbb{K}_X$-linear map of degree $1$ satisfying $d^2=0$. Similarly, for a chain complex, instead of a co-boundary map we have a boundary map $\partial: \mathcal{C}^\bullet \rightarrow \mathcal{C}^{\bullet -1}$, a $\mathbb{K}_X$-linear map of degree $-1$ satisfying $\partial^2=0$.
To address the following topics of discussion consistently, we consider such (co)chain complexes where $\mathcal{O}$-linearity of the (co)boundary maps is not necessarily required.

For a $\mathcal{O}$-module $\mathcal{C}$, we can canonically form the following sheaves of graded vector spaces $\wedge^\bullet_{\mathcal{O}}$ $\mathcal{C}:=\oplus_{i} \wedge^i_{\mathcal{O}} \mathcal{C}$ and $\mathcal{C}^{\otimes^ \bullet}:= \oplus_i$ $\mathcal{C}^{\otimes ^i}$, where $i \geq 0$ or $i \in \mathbb{N} \cup \{0\}$. These form exterior algebra and tensor algebra  of $\mathcal{C}$ with respect to the products $\wedge$ and $\otimes$ over $\mathcal{O}$ respectively.

	For an $\mathcal{O}$-module $\mathcal{F}$ on $X$, there exists an injective resolution  $\mathcal{F} \overset{\sim}{\rightarrow} \mathcal{C}^\bullet(\mathcal{F})$ $($a quasi-isomorphism) of $\mathcal{O}$-modules, known as the flabby Godement resolution of $\mathcal{F}$. Thus, for a complex of $\mathcal{O}$-modules $\mathcal{F}^\bullet$, we consider the associated bicomplex of sheaves $\mathcal{C}^\bullet(\mathcal{F}^\bullet)=(\mathcal{C}^p(\mathcal{F}^q))$ $(p,q \in \mathbb{N} \cup \{0\})$ of injective resolutions. The original complex is embedded in the total complex $\mathcal{K}^\bullet= tot(\mathcal{C}^\bullet(\mathcal{F}^\bullet))$, and this embedding is a quasi-isomorphism. The cohomology of the associated complex $\mathcal{K}^\bullet(X)=tot(\mathcal{C}^\bullet(\mathcal{F}^\bullet))(X)$ of global sections is called the hypercohomology of $\mathcal{F}^\bullet$, and denoted by $\mathbb{H}^\bullet(X, \mathcal{F}^\bullet)$ (see \cite{VID}).

	Denote the category of cochain complexes  of sheaves of $\mathcal{O}_X$-modules on an $a$-space $(X,\mathcal{O}_X)$ by $Shv(X)^\bullet$ and the category of $\mathbb{K}$-vector spaces by $Vect$.
	For $\mathcal{F}^\bullet\in Shv(X)^\bullet$, the $k$-th cohomology sheaf of the complex $\mathcal{F}^\bullet$ is
	\begin{center}
		$\mathscr{H}^k(\mathcal{F}^\bullet) := \mathscr{K}er (\mathcal{F}^k \rightarrow \mathcal{F}^{k+1})/{\mathscr{I}m(\mathcal{F}^{k-1} \rightarrow \mathcal{F}^{k})}$
	\end{center}
	$($where $\mathcal{F}^{-1}:=0)$ for $k \in \mathbb{N} \cup \{0\}$, considers as quotient of $\mathbb{K}_X$-vector spaces.
	Moreover, a map of complexes $\mathcal{F}^\bullet \rightarrow \mathcal{G}^\bullet$ is a quasi-isomorphism if the induced map on cohomology sheaves $\mathscr{H}^k(\mathcal{F}^\bullet) \rightarrow \mathscr{H}^k(\mathcal{G}^\bullet)$ is an isomorphism for all $k$.

	The $k$-th hypercohomology is a functor $\mathbb{H}^k(X,-):Shv(X)^\bullet \rightarrow Vect$ for $k \in \mathbb{N} \cup \{0\}$ (see \cite{VID,MS} for details) satisfying the following two conditions: \\
	A quasi-isomorphism of complexes $f^\bullet: \mathcal{F}^\bullet \rightarrow \mathcal{G}^\bullet$ induces an isomorphism $\mathbb{H}^k(X, f^\bullet)$, and 
	if $\mathcal{I}^\bullet$ is a complex of injective sheaves then  $\mathbb{H}^k(X, \mathcal{I}^\bullet) = H^k(\Gamma(X, \mathcal{I}^\bullet))$, where $\Gamma(X, \mathcal{I}^\bullet)$ is the space of global sections of $\mathcal{I}^\bullet$.

	Denote the graded vector space $\bigoplus\limits_{n=0}^{\infty} \mathbb{H}^n(X,\mathcal{F}^\bullet)$ as $\mathbb{H}^\bullet(X,\mathcal{F}^\bullet)$.
\begin{Rem}
	We use sheaf-theoretic generalizations of the well-known derived functors $Ext$ and $Cotor$ $($see \cite{GR},\cite{KP}$)$, by sheafifying the standard cochain complexes,  and subsequently considering the associated hypercohomologies.
\end{Rem}
\subsection{Chevalley-Eilenberg-de Rham complex and Koszul-Rinehart resolution}
For a Lie algebroid $\mathcal{L}$ over $(X, \mathcal{O}_X)$ with a representation $(\mathcal{E}, \nabla)$, we construct a cochain complex consisting of $\mathcal{O}_X$-modules equipped with a $\mathbb{K}_X$-linear degree $1$ differential. This serves as a generalization of the well-known Chevalley-Eilenberg-de Rham complex \cite{KM}. On the other hand, we consider a canonical resolution of $\mathcal{O}_X$ in the category of left $\mathscr{U}(\mathcal{O}_X, \mathcal{L})$-modules. Extending the standard correspondence established by Rinehart \cite{GR} for Lie-Rinehart algebras, we adapt it to the algebro-geometric setting.

 For a Lie algebroid $\mathcal{L}:=(\mathcal{L}, [\cdot, \cdot], \mathfrak{a})$ over $(X,\mathcal{O}_X)$, first recall some of the key notions in the following.

 Denote $\Omega^k_{\mathcal{L}(U)}(\mathcal{E}(U))$ or $^p Hom_{\mathcal{O}_X(U)}(\wedge^\bullet_{\mathcal{O}_X(U)}\mathcal{L}(U), \mathcal{E}(U))$ by the 
   module of $\mathcal{O}_X(U)$-linear maps from $\wedge^\bullet_{\mathcal{O}_X(U)}\mathcal{L}(U)$ to $\mathcal{E}(U)$ which induces $\mathcal{O}_X(V)$-module homomorphisms from $\wedge^\bullet_{\mathcal{O}_X(V)}\mathcal{L}(V)$ to $\mathcal{E}(V)$ for each open subset $V$ in $U$,
   and are compatible with restriction maps of the presheaf $U \mapsto \wedge^\bullet_{\mathcal{O}_X(U)}\mathcal{L}(U)$ and of the sheaf $\mathcal{E}$.
   Thus, in the underlying process, we consider the presheaf $^p \wedge^{\bullet}_{\mathcal{O}_X} \mathcal{L}$ given by $U \mapsto \wedge^{\bullet}_{\mathcal{O}_X(U)} \mathcal{L}(U)$ together with all the presheaf homomorphisms from $^p \wedge^{\bullet}_{\mathcal{O}_X} \mathcal{L}$ to $\mathcal{E}$. Consider the following cochain complex
\begin{align}\label{Chevalley-Eilenberg-de Rham complex} 
	\Omega^{\bullet}_{\mathcal{L}}(\mathcal{E}):= (\mathscr{H}om_{\mathcal{O}_X}(\wedge^\bullet_{\mathcal{O}_X}\mathcal{L}, \mathcal{E}), d).
\end{align}
\label{pre CED}
It is the sheafification of the presheaf of cochain complexes (see \cite{GR,CW3})
$$U\mapsto (^p Hom_{\mathcal{O}_X(U)}(\wedge^\bullet_{\mathcal{O}_X(U)}\mathcal{L}(U), \mathcal{E}(U)), d_U)=:\Omega_{\mathcal{L}(U)}^\bullet(\mathcal{E}(U)),$$
where the differential $d_U$ associated with an open subset $U$ of $X$ is given by
\begin{align*}
	\begin{split}
		d_U(\omega)(D_1\wedge \cdots \wedge D_{k+1}) & = \sum^{k+1}_{i=1}(-1)^{i+1}~\nabla^U_{D_i}(\omega(D_1\wedge \cdots \wedge \hat{D_i} \wedge \cdots \wedge D_{k+1}))\\
		&+ \sum_{i< j}(-1)^{i+j}~\omega([D_i, D_j]\wedge D_1\wedge \cdots \wedge \hat{D_i}\wedge \cdots \wedge \hat{D_j}\wedge \cdots \wedge D_{k+1}),
	\end{split}
\end{align*}
where $D_1, \dots, D_{k+1} \in \mathcal{L}(U)$ and
$  \omega \in \Omega^k_{\mathcal{L}(U)}(\mathcal{E}(U))$, 
and $\nabla:~ U \mapsto \nabla^U$ (see the notion in (\ref{L-connection})) is the sheafification.  Notice that, the differential $d$ (i.e. $d^2=0)$ is a $\mathbb{K}_X$-linear map, but not an $\mathcal{O}_X$-module homomorphism $($satisfies the Leibniz rule: $d(f~ \omega)= df \wedge \omega + f~ d\omega$, for $f \in \mathcal{O}_X(U),~ \omega \in \Omega^k_{\mathcal{L}(U)}(\mathcal{E}(U)))$.

In addition, if $\mathcal{L}$ is a locally free $\mathcal{O}_X$-modules of finite rank, then we get
\begin{align}\label{special case}
	\mathscr{H}om_{\mathcal{O}_X}(\wedge^{\bullet}_{\mathcal{O}_X}\mathcal{L}, \mathcal{E}) \cong \wedge^{\bullet}_{\mathcal{O}_X}\mathcal{L}^* \otimes_{\mathcal{O}_X} \mathcal{E},
\end{align}
which implies $\Omega^{\bullet}_{\mathcal{L}}(\mathcal{E}) \cong \Omega^{\bullet}_{\mathcal{L}}\otimes_{\mathcal{O}_X} \mathcal{E}$ where $\Omega^{\bullet}_{\mathcal{L}}:=\Omega^{\bullet}_{\mathcal{L}}(\mathcal{O}_X)$. This complex is equivalent to the standard complex associated with a smooth Lie algebroid \cite{KM} or a holomorphic Lie algebroid \cite{BP}. In particular, we get the Chevalley-Eilenberg complex when $\mathcal{L}$ is a Lie algebroid over a point; and the de Rham complex $\Omega^\bullet_X$ when  $\mathcal{L}=\mathcal{T}_X$ and $\mathcal{E}=\mathcal{O}_X$ over a non-singular (or smooth) $a$-space $X$ $($see Remark \ref{special $a$-spaces}$)$.

The associated hypercohomology of the cochain complex (\ref{Chevalley-Eilenberg-de Rham complex}) of $\mathcal{O}_X$-modules is called the Lie algebroid cohomology of $\mathcal{L}$ with coefficient in $\mathcal{E}$ and denoted by $\mathbb{H}^\bullet(\mathcal{L}, \mathcal{E})$.

Consider the chain complex of left $\mathscr{U}(\mathcal{O}_X, \mathcal{L})$-modules, called the Koszul-Rinehart complex: 
\begin{align} \label{Koszul-Rinehart complex}
	\mathscr{K}_{\mathcal{O}_X}^\bullet \mathcal{L}:=  (\mathscr{U}(\mathcal{O}_X, \mathcal{L})\otimes_{\mathcal{O}_X} \wedge^\bullet_{\mathcal{O}_X} \mathcal{L}, \partial),
\end{align} it is the sheafification of the presheaf of cochain complexes $($see \cite{GR,CW3}$)$ \label{pre KR}
$$U \mapsto (\mathcal{U}(\mathcal{O}_X(U), \mathcal{L}(U))\otimes_{\mathcal{O}_X(U)} \wedge^\bullet_{\mathcal{O}_X(U)} \mathcal{L}(U), \partial_U)=:K^\bullet_{\mathcal{O}_X(U)}\mathcal{L}(U)$$ where the differential $\partial_U$ associated with an open subset $U$ of $X$ is given by
\begin{align*} 
	\begin{split}
		\partial_U(\tilde{D}\otimes D_1\wedge \cdots \wedge D_{k}) & = \sum^k_{i=1}(-1)^{i-1}\tilde{D} \bar{D_i}\otimes D_1\wedge \cdots \wedge \hat{D_i} \wedge \cdots D_k\\
		& +  \sum\limits_{i< j}(-1)^{i+j}\tilde{D}\otimes [D_i, D_j]\wedge D_1\wedge \cdots \wedge \hat{D_i} \wedge \cdots \wedge \hat{D_j}\wedge \cdots \wedge D_{k},
	\end{split}
\end{align*}
where $\tilde{D} \bar{D_i}$ define by the canonical map $\iota_U:\mathcal{L}(U) \rightarrow \mathcal{U}(\mathcal{O}_X(U),\mathcal{L}(U))$, for all $\tilde{D}\in \mathcal{U}(\mathcal{O}_X(U),\mathcal{L}(U))$ and $D_1, \dots, D_k \in \mathcal{L}(U)$. Notice that the boundary (or chain) map $\partial$ is $\mathscr{U}(\mathcal{O}_X, \mathcal{L})$-linear.

When $\mathcal{L}$ is locally free as $\mathcal{O}_X$-module, it provides a locally free resolution of $\mathcal{O}_X$ as a left $\mathscr{U}(\mathcal{O}_X, \mathcal{L})$-module (implies the stalk-wise exactness), known as the Koszul-Rinehart resolution. That is, it is given by the (augmented) chain complex of left $\mathscr{U}(\mathcal{O}_X, \mathcal{L})$-modules 
\begin{align}\label{Koszul-Rinehart}
	\mathscr{K}_{\mathcal{O}_X}^\bullet \mathcal{L} \rightarrow \mathcal{O}_X,
\end{align}
defined by the counit map (or augmentation map) $\epsilon : \mathscr{U}(\mathcal{O}_X, \mathcal{L}) \rightarrow \mathcal{O}_X$ (is given as $\epsilon(\tilde{D})=\tilde{D}(1)$ for $\tilde{D}\in  \mathscr{U}(\mathcal{O}_X, \mathcal{L})$), provides the quasi-isomorphism.

\begin{Rem} \label{dual K-R}
	In the last section we consider analogous resolutions for $\mathcal{O}_X$ in the category of left $\mathcal{S}_{\mathcal{O}_X} \mathcal{L}$-comodules, where we replace the complex $\mathscr{K}_{\mathcal{O}_X}^\bullet \mathcal{L}$ with its dual complex $($see (\ref{Koszul 1})$)$, the dual Koszul-Rinehart complex $\widetilde{\mathscr{K}}_{\mathcal{O}_X}^\bullet\mathcal{L}$ $($sheafifying the local counterpart as described in \cite{KP} accordingly$)$.
\end{Rem}
\begin{Rem}
	In \cite{BP,BRT,UB}, the Lie algebroid cohomology  have considered for a special kind of Lie algebroid $\mathcal{L}$ (consists by coherent sheaf or locally free sheaf of finite rank) over a Noetherian separated scheme or a complex manifold $(X, \mathcal{O}_X)$. The Lie algebroid cohomology of $\mathcal{L}$ with values in an $\mathcal{L}$-module $\mathcal{E}$ (a coherent $\mathcal{O}_X$-module) is the hypercohomology of the complex $(\wedge^\bullet_{\mathcal{O}_X}\mathcal{L}^* \otimes_{\mathcal{O}_X} \mathcal{E},$ $d)$. In particular, this appears as a special case (\ref{special case}) of the complex described in (\ref{Chevalley-Eilenberg-de Rham complex}).
\end{Rem}
\begin{Rem}
Moreover, one can consider the cochain complex of K\"ahler differentials 
\((\wedge^{\bullet}_{\mathcal{O}_X} \Omega^1_X, d_X)\) over a scheme \((X, \mathcal{O}_X)\), 
where \(\Omega^1_X\) is the \(\mathcal{O}_X\)-module of K\"ahler differential 1-forms on \(X\), and the differential \(d_X\) extends the canonical map \(d^1_X: \mathcal{O}_X \to \Omega^1_X\). 
If \(\Omega^1_X\) is reflexive, then \((\mathcal{T}_X)^* = (\Omega^1_X)^{**} \cong \Omega^1_X\) as \(\mathcal{O}_X\)-modules, a property that holds when $X$ is a smooth scheme of finite type. Consequently, the algebraic de Rham complex of $X$ is recovered by substituting \(\mathcal{L} = \mathcal{T}_X\) in the Lie algebroid complex discussed in \cite{UB}.
\end{Rem}
\begin{Nota} \label{Special open} From now on we use a special type of open cover $\{U_x \mid  x\in X \}$ of $X$ for a locally free Lie algebroid $\mathcal{L}$, where the restrictions $\mathcal{L}|_{U_x}$ are free $\mathcal{O}_X|_{U_x}$-modules for every $U_x$. These open sets $U_x$'s are said to be special open sets.
\end{Nota}

We now describe Lie algebroid cohomology in terms of a derived functor. This result is more general than the corresponding result in \cite{UB, BP}, where underlying $\mathcal{O}_X$-module of a Lie algebroid is assumed to be a coherent sheaf. In contrast, we consider the more general case in which it is only assumed to be quasi-coherent. Moreover, our proof differs from that given in \cite{UB}.

\begin{Thm} \label{Lie algebroid Cohomology as derived functor}
	Let $\mathcal{L}$ be a locally free Lie algebroid over $(X, \mathcal{O}_X)$ and $\mathcal{E}$ a representation of $\mathcal{L}$. Then we get an isomorphism of graded vector spaces
	\begin{center}
		$\mathbb{H}^\bullet (\mathcal{L}, \mathcal{E})\cong Ext^\bullet_{\mathscr{U}(\mathcal{O}_X, \mathcal{L})}(\mathcal{O}_X, \mathcal{E})$. 
	\end{center}
\end{Thm}
\begin{proof}
	For every $x\in X$, the space of sections $\mathcal{L}(U_x)$ is a $(\mathbb{K},\mathcal{O}_{X}(U_x))$-Lie-Rinehart algebra and a free module over $\mathcal{O}_{X}(U_x)$ (see Notation \ref{Special open}). Thus, the Koszul-Rinehart resolution (\ref{Koszul-Rinehart}) is a locally free resolution of the sheaf of local rings $\mathcal{O}_{X}$. Note that, here the $\mathcal{L}$-module $\mathcal{E}$ is also considered as a $\mathscr{U}(\mathcal{O}_X, \mathcal{L})$-module.

	Here, we consider the two naturally associated presheaf of cochain complexes $($see Section \ref{pre CED}, \ref{pre KR}$)$
	\begin{center}
		$U\mapsto \Omega_{\mathcal{L}(U)}^\bullet(\mathcal{E}(U))$,\\
		\vspace{.2 cm}
		\hspace{1 cm}	$U\mapsto Hom_{\mathcal{U}(\mathcal{O}_X(U),\mathcal{L}(U))}(K^\bullet_{\mathcal{O}_X(U)}\mathcal{L}(U), \mathcal{E}(U))$.
	\end{center}
	We have a canonical morphism between these two presheaves of cochain complexes as $\mathbb{K}_X$-modules, which is a local isomorphism or stalk-wise isomorphism, using results from \cite{GR} as local descriptions (applying on each special open neighborhoods $U_x$'s of $x \in X$). Thus, the associated sheaves of cochain complexes are isomorphic as $\mathbb{K}_X$-modules (using Lemma $3.2.$ and Lemma $3.3.$ from \cite{AA}), i.e.
	$$\Omega^\bullet_{\mathcal{L}} (\mathcal{E}) \cong \mathscr{H}om_{\mathscr{U}(\mathcal{O}_X, \mathcal{L})}(\mathscr{K}_{\mathcal{O}_X}^\bullet \mathcal{L}, \mathcal{E}).$$
	Hence, the Chevalley-Eilenberg-de Rham complex of $\mathcal{L}$ with coefficient in $\mathcal{E}$ is isomorphic to dual of the Koszul-Rinehart resolution of $\mathcal{O}_X$ by $\mathcal{E}$ in the category of left $\mathscr{U}(\mathcal{O}_X, \mathcal{L})$-modules.
	
	Therefore, the associated cohomology sheaves (see Section \ref{Hypercohomology}) are isomorphic, i.e.
	\begin{center}
		$\mathscr{H}^k(\Omega^\bullet_{\mathcal{L}}, \mathcal{E})\cong \mathscr{E}xt^k_{\mathscr{U}(\mathcal{O}_X, \mathcal{L})}(\mathcal{O}_X, \mathcal{E})$, for every $k$.
	\end{center}
	
	Now, applying the hypercohomology functor $\mathbb{H}^\bullet(X,-)$ (see Section  \ref{Hypercohomology}) we get the required result (since $\mathbb{H}^\bullet(\mathcal{L}, \mathcal{E})= \mathbb{H}^\bullet(X,\Omega^\bullet_{\mathcal{L}}(\mathcal{E}))$ and $Ext^\bullet_{\mathscr{U}(\mathcal{O}_X, \mathcal{L})}(\mathcal{O}_X, \mathcal{E})=\mathbb{H}^\bullet(X, \mathscr{H}om_{\mathscr{U}(\mathcal{O}_X, \mathcal{L})}(\mathscr{K}_{\mathcal{O}_X}^\bullet \mathcal{L}, \mathcal{E}))$).
\end{proof}
\begin{Rem} One can find a proof for locally free Lie algebroids $\mathcal{L}$ of finite rank over a Noetherian separated scheme $(X,\mathcal{O}_X)$ discussed in \cite{UB} by U. Bruzzo, by using the ideas of \cite{GR} on the level of stalks $\mathcal{L}_x$ $($since $\mathcal{L}_x$ is a $(\mathbb{K},\mathcal{O}_{X,x})$-Lie-Rinehart algebra, projective (in fact, free) module over $\mathcal{O}_{X,x}$ for all $x \in X)$. 
\end{Rem}

\begin{Cor}
	Let $(X, \mathcal{O}_X)$ be a complex manifold or smooth algebraic variety and $Y$ a free divisor in $X$. Therefore, the Lie algebroid $\mathcal{T}_X(-log Y)$ of logarithmic derivation is locally free $\mathcal{O}_X$-module (of finite rank). Using Example \ref{sheaf of log diff operators}, the logarithmic de Rham cohomology (hypercohomology of the complex (\ref{log de Rham complex}) described in Remark \ref{log de Rham cohomology 1} and Remark \ref{log de Rham cohomology 2} for some special cases, as given in (\ref{holomorphic Lie algebroid cohomology})) can be expressed as
	$$\mathbb{H}^\bullet(X, \Omega^\bullet_X(log Y)) \cong Ext^\bullet_{\mathcal{D}_X(-log Y)}(\mathcal{O}_X,\mathcal{O}_X).$$
\end{Cor}

\begin{Cor}
	Consider the path algebroid $\mathcal{P}_X$ for a non-singular special $a$-space $(X,\mathcal{O}_X)$  $($see Remark \ref{special $a$-spaces}$)$ and its universal enveloping algebroid $\mathbb{D}_X$ (see Example \ref{free Lie algebroid}). Note that, a module over the free Lie algebroid $\mathcal{P}_X$ can be viewed as a vector bundle with $\mathcal{T}_X$-connection $($not necessarily flat$)$, and a $\mathbb{D}_X$-module describes some system of linear PDEs on the space of paths (see \cite{MK} for details).
	
	The Theorem \ref{Lie algebroid Cohomology as derived functor} for the Lie algebroid $\mathcal{P}_X$ provides the isomorphism
	$$\mathbb{H}^\bullet (\mathcal{P}_X, \mathcal{E})\cong Ext^\bullet_{\mathbb{D}_X}(\mathcal{O}_X, \mathcal{E}).$$ 
\end{Cor}

\begin{Cor}
	Let $\mathcal{A}$ be a cocomplete locally graded free $\mathcal{O}_X/{\mathbb{K}_X}$-bialgebra \cite{AA} of finite or infinite type (for examples consider Remark \ref{cocomplete locally graded free bialgebras}). Then the Lie algebroid $\mathscr{P}(\mathcal{A})$ (sheaf of primitive elements, see Remark \ref{primitives}) is a locally free $\mathcal{O}_X$-module (of finite or infinite rank accordingly). Thus, we have $\mathcal{A}\cong \mathscr{U}(\mathcal{O}_X,\mathscr{P}(\mathcal{A}))$ \cite{AA}. Therefore, applying Theorem \ref{Lie algebroid Cohomology as derived functor} for $\mathcal{L}=\mathscr{P}(\mathcal{A})$ provides the isomorphism
	$$\mathbb{H}^\bullet (\mathscr{P}(\mathcal{A}), \mathcal{E})\cong Ext^\bullet_{\mathcal{A}}(\mathcal{O}_X, \mathcal{E}).$$
\end{Cor}
\subsubsection{On \v{C}ech cohomology} \label{Cech cohomology derived}
	In complex geometry $($algebraic geometry$)$, for a cochain complex of coherent sheaves $\mathcal{F}^\bullet$ over an analytic space $($algebraic variety$)$, we compute the hypercohomology $\mathbb{H}^\bullet(X,\mathcal{F}^\bullet)$ via \v{C}ech cohomology $\check{H}^\bullet(\mathcal{U},\mathcal{F}^\bullet)$ associated with some good open cover $\mathcal{U}=\{U_i\}_i$ of $X$ as (a canonical isomorphism) $$\mathbb{H}^\bullet(X,\mathcal{F}^\bullet)\cong \check{H}^\bullet(\mathcal{U},\mathcal{F}^\bullet).$$  Specifically, we provide a good open cover by connected Stein spaces (affine varieties) and we use Leray's theorem \cite{VID,SR} and Cartan-Serre's vanishing theorem \cite{FF} for sheaf cohomology.

	When we consider a locally free Lie algebroid $\mathcal{L}$ of finite rank over an analytic space (algebraic variety) $(X,\mathcal{O}_X)$, we can compute Lie algebroid cohomology as derived functor using \v{C}ech cohomology.
	
	In both of these cases we obtain the isomorphism
	$\mathbb{H}^\bullet (\mathcal{L}, \mathcal{E})\cong \check{H}^\bullet(\mathcal{U},\Omega^\bullet_{\mathcal{L}}(\mathcal{E}))$.
	
	Since over Stein spaces $($affine varieties$)$ $U_i$, the $\mathcal{O}_X(U_i)$-module $\mathcal{L}(U_i)$ is a projective module \cite{Mo, CW}, thus we get an isomorphism
	\begin{center}
		$H^\bullet(\mathcal{L}(U_i),\mathcal{E}(U_i)) \cong Ext^\bullet_{\mathcal{U}(\mathcal{O}_X(U_i),\mathcal{L}(U_i))}(\mathcal{O}_X(U_i),\mathcal{E}(U_i))$.
	\end{center}
	Hence, $$\check{H}^\bullet(\mathcal{U},\Omega^\bullet_{\mathcal{L}}(\mathcal{E})) \cong \check{H}^\bullet(\mathcal{U}, \mathscr{H}om_{\mathscr{U}(\mathcal{O}_X,\mathcal{L})}(\mathscr{K}_{\mathcal{O}_X}^\bullet \mathcal{L},\mathcal{E})).$$ 
	
	Also, note that the whole ideas works for the classical case of Lie algebroids over smooth manifolds (i.e. for smooth Lie algebroids). In \cite{BRT}, the Lie algebroid $($hyper$)$cohomology of holomorphic $($algebraic$)$ Lie algebroids over a complex manifold $($smooth scheme$)$ has expressed by \v{C}ech cohomology using a good open cover consisting with connected Stein manifolds $($affine schemes$)$. However, the case associated with analytic spaces is not considered. Additionally, this result can be extended to the quasicoherent $\mathcal{O}_X$-module cases.
	
	More generally, for a cochain complex of quasicoherent sheaves $\mathcal{F}^\bullet$ over a Noetherian separated scheme, by considering an affine open cover $\mathcal{U}=\{U_i\}_i$ of $X$, we get a canonical isomorphism $\mathbb{H}^\bullet(X,\mathcal{F}^\bullet)\cong \check{H}^\bullet(\mathcal{U},\mathcal{F}^\bullet)$ \cite{RH}. Thus, we get a similar result for a locally free Lie algebroid over a Noetherian separated scheme.

\subsection{Chern classes} \label{Chern classes}
Lie algebroid (or Lie-Rinehart algebra) connections, curvatures and the associated Characteristic classes (e.g. Picard group, Chern class) has been studied in the context of differential geometry \cite{RF} and algebric geometry \cite{JH,HM3,HM,HOM}. We review and extend some of the facts in our context as follows.

Classically, the first Chern class produces a bijection between the set of all isomorphism classes of complex line bundles on a manifold $X$ and the cohomology group $H^2(X, \mathbb{Z})$. More generally, for any Lie-Rinehart algebra $L$ which is projective as a left $R$-module there is a one-to-one correspondence between isomorphism classes of projective $R$-modules of finite rank and
the Lie-Rinehart cohomology group $H^2(L,R)$.

\subsubsection{Abelian Lie algebroid extensions}
Consider a $(\mathbb{K},R)$-Lie-Rinehart algebra $L$ which is  projective as an $R$-module. Let $(E,\nabla)$ be an $L$-module, then we have the isomorphism of $\mathbb{K}$-vector spaces  $$H^2(L,E)\cong Ext^1(L,E,\nabla),$$ where $Ext^1(L,E,\nabla)$ is the set of equivalence classes of extensions of $L$ by the flat connection $(E, \nabla)$ \cite{JH,HOM}.
This result can be viewed in the algebraic geometry set up (over the affine scheme $Spec R$).  We derive an analogous result for a locally free Lie algebroid $\mathcal{L}$ over an $a$-space $(X, \mathcal{O}_X)$ with an $\mathcal{L}$-module $(\mathcal{E}, \nabla)$ as follows.
We denote the set of equivalence classes of abelian extensions of a Lie algebroid $\mathcal{L}$ by a flat $\mathcal{L}$-connection $(\mathcal{E}, \nabla)$ as $Ext^1(\mathcal{L}, \mathcal{E},\nabla)$. An element $[\mathcal{L}'] \in Ext^1(\mathcal{L}, \mathcal{E},\nabla)$ represent a s.e.s. of Lie algebroids
$$\mathcal{E} \rightarrow \mathcal{L}'\rightarrow \mathcal{L}$$
upto isomorphism in the usual sense, where $\mathcal{E}$ is the trivial Lie algebroid. Thus, as a local description, on the level of stalks, for each point $x \in X$, we have the s.e.s. of $(\mathbb{K}, \mathcal{O}_{X,x})$-Lie-Rinehart algebras
$$\mathcal{E}_x \rightarrow {\mathcal{L}'}_x\rightarrow \mathcal{L}_x.$$
\begin{Thm} \label{extensions of L by E}
Let $(\mathcal{L}, [\cdot, \cdot], \mathfrak{a})$ be a locally free Lie algebroid  over an $a$-space $(X, \mathcal{O}_X)$ with an $\mathcal{L}$-module $(\mathcal{E}, \nabla)$. Then, we get a canonical isomorphism of $\mathbb{K}$-vector spaces
$$\mathbb{H}^2(\mathcal{L}, \mathcal{E})\cong Ext^1(\mathcal{L}, \mathcal{E},\nabla).$$
\end{Thm}
\begin{proof}
	Given a Lie algebroid $(\mathcal{L}, [\cdot, \cdot], \mathfrak{a})$ over $(X, \mathcal{O}_X)$ and an element $($an isorphism class$)$ $[\omega] \in \mathbb{H}^2(\mathcal{L}, \mathcal{E})$ where $\omega \in \Omega^2_{\mathcal{L}}(\mathcal{E})$, we get a Lie algebroid extension of $\mathcal{L}$ by $\mathcal{E}$ (upto isomorphism) as follows.
	
	Consider the $\mathcal{O}_X$-module $\mathcal{E} \oplus \mathcal{L}$ with the $\mathbb{K}_X$-Lie algebra structure given by
	$$[(s,D), (s',D')]_{\omega}:=(\nabla_D(s') - \nabla_{D'}(s) +\omega(D,D'), [D, D']),$$
	for $s, s'\in \mathcal{E}$ and $D, D' \in \mathcal{L}$, together with the map
	$$\mathfrak{a}_{\omega}:\mathcal{L}_{\omega} \rightarrow \mathcal{T}_X,~ (s, D) \mapsto \mathfrak{a}(D).$$
	Thus, $\mathcal{L}_{\omega}:=(\mathcal{E} \oplus \mathcal{L}, [\cdot, \cdot]_{\omega}, \mathfrak{a}_{\omega})$ forms a Lie algebroid over $(X, \mathcal{O}_X)$. 
     Let $[\omega]=[\omega']$ for some $\omega' \in \Omega^2_{\mathcal{L}}(\mathcal{E})$, i.e. $\omega - \omega'=d\theta$ for some $\theta \in \Omega^1_{\mathcal{L}}(\mathcal{E})$. Then the map 
	$\tilde{\theta}:\mathcal{L}_{\omega} \rightarrow \mathcal{L}_{\omega'}$ defined by $\tilde{\theta}(s+D)=(s+\theta(D))+D$ forms an isomorphism of Lie algebroids. Thus,
	we get an abelian Lie algebroid extension (as an equivalence class) of $\mathcal{L}$ by the module $(\mathcal{E}, \nabla)$ i.e. a s.e.s. of Lie algebroids where $\mathcal{E}$ is the trivial Lie algebroid
	$$\mathcal{E} \hookrightarrow \mathcal{L}_{\omega} \twoheadrightarrow \mathcal{L},$$
	and the resultant isomorphism class $[\mathcal{L}_{\omega}] \in Ext^1(\mathcal{L}, \mathcal{E},\nabla)$.

	Consider an abelian Lie algebroid extension $\mathcal{L}'$ of $\mathcal{L}$ by $\mathcal{E}$ upto isomorphism, i.e. $[\mathcal{L}'] \in Ext^1(\mathcal{L}, \mathcal{E},\nabla)$ as
	\begin{align} \label{extension}
	 \mathcal{E} \rightarrow \mathcal{L}'\rightarrow \mathcal{L},
	\end{align}
	where $(\mathcal{E}, \nabla)$ is an $\mathcal{L}$-module.
	Since $\mathcal{L}$ is locally free $\mathcal{O}_X$-module, the s.e.s. (\ref{extension}) splits as an $\mathcal{O}_X$-module, i.e. there is an isomorphism $\tilde{\phi}: \mathcal{L}' \rightarrow \mathcal{E} \oplus \mathcal{L}$ for the split $($or section map$)$ $\phi: \mathcal{L} \rightarrow \mathcal{L}'$. Thus, the possible choice for a $\mathbb{K}_X$-Lie bracket on $\mathcal{E} \oplus \mathcal{L}$ is given as follows:
	$$[(s,D), (s',D')]':=[s, s']'+ [D, s']'-[D', s]'+[D, D']',$$
	 where $[s, s']'=0, ~[D, s']'=\nabla_D(s'),~ [D', s]'=\nabla_{D'}(s)$, 
	 ~$[D, D']'=\omega(D, D')+ [D, D']$ for some $\omega \in \Omega^2_{\mathcal{L}}(\mathcal{E})$ with $d(\omega)=0$,~ for $s, s' \in \mathcal{E},~ D, D' \in \mathcal{L}$. Let $\phi': \mathcal{L} \rightarrow  \mathcal{L}'$ be another such splitting map  and 
	 $\omega' \in \Omega^2_{\mathcal{L}}(\mathcal{E})$ associated to the isomorphism  $\tilde{\phi'}: \mathcal{L}' \rightarrow \mathcal{E} \oplus \mathcal{L}$. Thus, using the automorphism $\tilde{\phi}' \circ \tilde{\phi}^{-1}$ on $\mathcal{E} \oplus \mathcal{L}$ we get $\omega - \omega'= d(\phi - {\phi}')$.
	 Thus, we get the required element $[\omega] \in \mathbb{H}^2(\mathcal{L}, \mathcal{E})$ associated to the equivalence class $[\mathcal{L}']$. 
	 
	 This one-to-one correspondence between the $\mathbb{K}$-vector spaces 
	$\mathbb{H}^2(\mathcal{L}, \mathcal{E})$ and $Ext^1(\mathcal{L},  \mathcal{E},\nabla)$ (with the standard vector space structure) can be extended to an isomorphism.
\end{proof}
  \begin{Rem}
  	In particular, for $\mathcal{E}= \mathcal{O}_X$ with the standard $\mathcal{L}$-module structure given by the anchor map, this kind of result is given in the holomorphic context using \v{C}ech cohomology \cite{PT}. 
  \end{Rem}
  \subsubsection{First Chern class} \label{first Chern class}
Note that the Chern class of a finitely generated projective $R$-module $E$ with an $L$-connection $\nabla$ on $E$, defined by the trace of the curvature $R_{\nabla}$ associated with the connection $\nabla$ \cite{HM,HOM}. If $R$ is a regular $\mathbb{K}$-algebra (implies $L:=Der_{\mathbb{K}}(R)$ is a finitely generated projective $R$-module) then by considering the  covariant connection on $E$ we get the first Chern class  given by trace of the curvature $$c_1(E) =Trace(R_{\nabla}) \in H^2(Der_{\mathbb{K}}(R),R).$$ 
Analogous result holds for a locally free $\mathcal{O}_X$-module $\mathcal{E}$ of finite rank over the affine scheme $X=Spec R$.
We compute analogous result in the general set up by sheafifying the local descriptions given in \cite{HM3,HM}.

For a complex manifold $X$, a locally free $\mathcal{O}_X$-module $\mathcal{E}$ of finite rank (equivalently a holomorphic vector bundle over $X$) always have covariant derivative as a flat  $\mathcal{T}_X$-connection on $\mathcal{E}$. 
An arbitrary $\mathcal{T}_X$-connection on $\mathcal{E}$ is given by an $\mathcal{O}_X$-linear map $$\nabla : \mathcal{T}_X \rightarrow \mathcal{E}nd_{\mathbb{C}_X}(\mathcal{E})$$ 
satisfying the Leibniz rule, and its curvature is the $\mathcal{O}_X$-linear map $$R_{\nabla}:\wedge^2_{\mathcal{O}_X} \mathcal{T}_X \rightarrow \mathcal{E}nd_{\mathbb{C}_X}(\mathcal{E})$$ defined as $$R_{\nabla}(D\wedge D')= [\nabla_D,\nabla_{D'}]_c - \nabla_{[D,D']_c}$$ for any two sections $D,D'$ of $\mathcal{T}_X$. In general, for an $\mathcal{L}$-connection $\nabla$ on $\mathcal{E}$, one verifies that the curvature $R_{\nabla} \in \mathscr{H}om_{\mathcal{O}_X}(\wedge^2_{\mathcal{O}_X}\mathcal{L}, \mathcal{E}nd_{\mathcal{O}_X}(\mathcal{E}))$, where $\mathcal{L}$ is a Lie algebroid.

Consider the induced $\mathcal{T}_X$-connection $ad~\nabla=[\nabla, \cdot]$ on $\mathcal{E}nd_{\mathcal{O}_X}(\mathcal{E})$, a flat connection.
Thus by considering the Chevalley-Eilenberg-de Rham complex with coefficient in $\mathcal{E}nd_{\mathcal{O}_X}(\mathcal{E})$, 
we get the cohomology class $[R_{\nabla}] \in \mathbb{H}^2(X, \mathcal{E}nd_{\mathcal{O}_X}(\mathcal{E}))$ (since $d^2(R_{\nabla})=0)$. Hence, the Chern class of the $\mathcal{T}_X$-module $(\mathcal{E}, \nabla)$ is $$c_1(\mathcal{E}):= Trace(R_{\nabla}) \in \mathbb{H}^2(X, \mathcal{O}_X)$$ where $Trace:\mathcal{E}nd_{\mathcal{O}_X}(\mathcal{E}) \rightarrow \mathcal{O}_X$ is the trace map.

The restriction of a $\mathcal{T}_X$-connection  to the Lie algebroid $\mathcal{L}:=\mathcal{T}_X(- log~ Y)$ $($i.e. $\Omega_{\mathcal{L}}^1 := \Omega_X^1( log~ Y)$) induces a logarithmic connection which is a meromorphic connection with simple poles along the divisor $Y$. Moreover, if $Y$ is a free divisor, then for a locally free $\mathcal{O}_X$-module $\mathcal{E}$ of finite rank with an $\mathcal{L}$-connection we get the first Chern class $c_1(\mathcal{E}) \in \mathbb{H}^2(\mathcal{L}, \mathcal{O}_X)$. 

We can extend the notion of Chern classes for $\mathcal{T}_X$ or $\mathcal{T}_X(- log~ Y)$ to arbitrary Lie algebroids. Note that,  an $\mathcal{L}$-connection does not always exist in general, but it does exist for a large class of cases \cite{PT, BRT, HOM}.

\section{Logarithmic de Rham cohomology } \label{Sec 4} 
We compute algebraic (analytic) de Rham cohomology groups of a family of non-singular varieties $Y_t$ for $t \in \mathbb{C}\setminus \{0\}$ and their associated singular variety $Y_0$ (which appears as principal divisors \cite{AA}). Thus, we derive cohomology of the Lie algebroid $\mathcal{T}_{Y_t}$ for $t \in \mathbb{C}\setminus\{0\}$ and study the singular case analogously.

Then, we compute hypercohomology of the logarithmic de Rham complexes \cite{Mor}, which simplifies to the Lie algebroid cohomology for $\mathcal{T}_X(- log Y_t)$ \cite{BP}, as we work over some free divisors $Y_t \subset X$ for $t \in \mathbb{C}$ \cite{CMD,BP}.

We use the notion $\langle \{s_1, \dots, s_n \} \rangle$ for $R$-module generated by $s_1, \dots, s_n \in E$, where $R$ is a $\mathbb{K}$-algebra and $E$ is an $R$-module (also extend to sheaf theoretic settings).

$(1)$ \textbf{Rectangular Hyperbolas:}\label{hyperbola} We consider the affine space $X:=\mathbb{A}^2$ ($\mathbb{C}^2$ with the Zariski topology) with its coordinate ring $\tilde{R}=\mathcal{O}_X(X)=\mathbb{C}[x,y]$. The principal ideal $I=\langle xy-t \rangle \subset \tilde{R}$ provides rectangular hyperbolas $Y_t:=V(I)$ (vanishing set of $I$) for every $t$ in $\mathbb{C}\setminus \{0\}$. For $t \in \mathbb{C} \setminus \{0\}$,  $Y_t$ is a non singular affine variety with its coordinate ring $R_t:=\mathcal{O}_{Y_t}(Y_t)=\mathbb{C}[x,y]/ \langle xy-t \rangle$. 

Here, $Y:=Y_1$ is homeomorphic to $\mathbb{C} \setminus \{0\}$ and $\mathbb{C} \setminus \{0\}$ is  same homotopy type of $S^1$. Thus, on considering the singular cohomology we get that 
\begin{center}
	$H^i(Y, \mathbb{C})=H^i(S^1, \mathbb{C})=  \mathbb{C} $ for $i=0,1$\\
	\hspace{-.6 cm}and $H^i(Y, \mathbb{C})=H^i(S^1, \mathbb{C})= \{0\}$ for $i \geq 2$.
\end{center}
We can use algebraic de Rham theorem \cite{VID,MS} to get algebraic de Rham cohomology of $Y$ by the singular cohomology of $S^1$. Now we recall the computation which helps us to follow the associated singular case when $t=0$ (which is a normal crossing divisor). 

On $Y$, the differential $df=xdy+ydx=0$ (where $f=xy-1$) and $\frac{1}{x} \in R_1= \mathbb{C}[x, x^{-1}]=:R$, hence $dy=- \frac{dx}{x^2}\in \langle dx \rangle$. Thus the space of differential $1$-forms on $Y$ (similar to K$\ddot{a}$hler differentials for $Spec~R$) is 
$$\Omega^1_R:=\frac{\langle \{dx,dy\} \rangle}{\langle df \rangle}=\langle dx \rangle.$$ Note that $\frac{dx}{x}|_Y$ is an algebraic differential $1$-form on $Y$ (here $\Omega^1_R$ is a free $R$-module of rank $1$). The space of algebraic differential $2$-forms on $Y$ is $\Omega^2_R:=\wedge^2_R \Omega^1_R$.  Now, $dx \wedge dy=dx \wedge - \frac{dx}{x^2}=0$ on $Y$ (since $y=\frac{1}{x}$ in $Y$) and thus we get $\Omega^2_R=\{0\}$.

Thus,  the algebraic de Rham complex of $Y$ (or $Spec R$) is  
$$R  \overset{d^0}{\longrightarrow} \Omega^1_R \overset{d^1}{\longrightarrow} 0  \cdots ,$$
where the first differential is given by $x\mapsto dx$ and $\frac{1}{x}\mapsto - \frac{dx}{x^2}$. Here, $\frac{dx}{x}$ is in the kernel of $d^1$ but $\frac{dx}{x}$ is not in the image of  $d^0$ 
(the only possible choice of preimage is $log x$, but it is not a polynomial or regular function). It follows that $H^1_{dR}(Y)$ is a $\mathbb{C}$-vector space of dimension $1$. Now, $H_{dR}^0(Y) \cong \mathbb{C}$ and $H_{dR}^n(Y)=\{0\}$ for $n \geq 2$. 

This is studied in \cite{AH}, as a topological invariant for non singular spaces, but the analogous computation for the associated singular case is not derived there.

Note that the affine algebraic set $Y$ can be considered as an analytic space. Moreover, this space can be viewed as a principal divisor \cite{AA,BP}.

We consider the associated logarithmic de Rham cohomology \cite{CMD,Mor} in the context of complex geometry.

\begin{Exm} \label{log de Rham cohomology 1} 
 \normalfont{Consider the sheaf of holomorphic functions $\mathcal{O}_X$  on $X:=\mathbb{C}^2$ and the ideal sheaf $\mathcal{I}=\langle xy-1 \rangle$ of $\mathcal{O}_X$. The analytic subspace  associated with the principal ideal sheaf $\mathcal{I}$  is $(Y:=V(\mathcal{I}),~ \mathcal{O}_Y:= \mathcal{O}_X/{\mathcal{I}})$ of $(X, \mathcal{O}_X)$.
	The sheaf of logarithmic derivations for $Y$ in $X$ is $$\mathcal{T}_X(-log~ Y)=\langle \{x\partial_x-y\partial_y, f ~\nabla f\}  \rangle,$$ where $\nabla f=y \partial_x + x \partial_y$ $($since $f=xy-1)$. The sheaf of logarithmic $1$-forms for $Y$ in $X$ (i.e. meromorphic $1$-forms $\omega$ of $X$ with poles along the divisor $Y$  such that $f \omega \in \Omega^1_X$ and $df \wedge \omega \in \Omega^2_X$ \cite{Mor}) is 
	$$\Omega^1_X(log~ Y)=\mathscr{H}om_{\mathcal{O}_X}(\mathcal{T}_X(-log~ Y),\mathcal{O}_X)=\langle \{dx,dy,\frac{df}{f}\} \rangle,$$ where $df=y dx+x dy$. 
	Note that, the relation $y~dx+x~dy-(xy-1)~\frac{df}{f}=0$ holds.

	The sheaf of logarithmic differential $2$-forms on $X$ associated with $Y$ is $$\Omega^2_X(log~ Y):=\wedge^2_{\mathcal{O}_X} \Omega^1_X(log~ Y)=\langle \{dx \wedge \frac{df}{f},dy \wedge \frac{df}{f},dx \wedge dy\} \rangle.$$ 
	
	To determine hypercohomology \cite{VID,MS} of the logarithmic de Rham complex
	\begin{center}
		$\Omega^{\bullet}_X(log~ Y):$	$\mathcal{O}_X \overset{d^0}{\longrightarrow}\Omega^1_X(log~ Y)\overset{d^1}{\longrightarrow}\wedge^2_{\mathcal{O}_X}\Omega^1_X(log~ Y) \overset{d^2}{\longrightarrow}\cdots$,
	\end{center}
	we have $\mathscr{K}er$ $d^0=\mathbb{C}_X$; $\mathscr{K}er$ $d^1$  is the $\mathbb{C}_X$-vector space generated by 
	\begin{center}
		$\{dg~| ~g \in \mathcal{O}_X\}\cup \{\frac{df}{f}\}$, i.e. $\mathscr{K}er$ $d^1=\mathscr{I}m$ $d^0 \oplus \mathbb{C}_X\{\frac{df}{f}\}$ (here $d^0=d)$;
	\end{center}
	and $\mathscr{K}er$ $d^2$ is the $\mathbb{C}_X$-vector space generated by $\{dg \wedge \frac{df}{f}~| ~g \in \mathcal{O}_X \} \cup \{dx \wedge dy\}= \Omega^2_X(log~ Y)$. Therefore,
	\begin{center}
		$ \mathscr{K}er$ $d^2= \mathbb{C}_X \{d^1(g~ \frac{df}{f})~|~ g \in \mathcal{O}_X\}\oplus \mathbb{C}_X  \{d^1(x~dy)\}=\mathscr{I}m$ $d^1$.
	\end{center} 
	
	Therefore, the hyperchomology groups (taking the space of global sections of the complex, since $X$ is a Stein manifold) are given by the following
	$$\mathbb{H}^0(X,\Omega^{\bullet}_X(log~ Y))=\mathbb{C}, ~ \mathbb{H}^1(X,\Omega^{\bullet}_X(log~ Y))=\mathbb{C}, ~ \mathbb{H}^n(X,\Omega^{\bullet}_X(log~ Y))=\{0\} ~(n \geq 2).$$

	On $Y$ we have $df=y dx+x dy=0$  implies $\frac{dx}{x}+\frac{dy}{y}=0$ ($x\neq 0, y\neq 0$) and by integrating it we get $xy=t$ are the solutions which parametrized by $t\in \mathbb{C}\setminus \{0\}$.
	
	Thus, we can view the normal crossing divisor $Y_0 :=\{(x,y)\in \mathbb{C}^2 \mid xy=0\}$ (a singular analytic space and a free divisor) as deformation of family of rectangular hyperbolas $Y_t$ (appears in \cite{MP} as deformation of a scheme). }
	
\end{Exm}
$(2)$ \textbf{Normal Crossing Divisor:} Here we compute the algebraic de Rham cohomology of the normal crossing divisor $Y_0$ which is a singular affine algebraic set in $X$ with the space of global section of its structure sheaf 
$$R_0 :=\mathcal{O}_{Y_0}(Y_0)=\frac{\mathbb{C}[x,y]}{\langle xy \rangle}.$$

The algebraic de Rham complex for $Y_0$ (or $Spec$ $R_0)$ is  
$$R_0  \overset{d}{\longrightarrow} \Omega^1_{R_0} \overset{d^1}{\longrightarrow} \Omega^2_{R_0} \overset{d^2}{\longrightarrow} 0 \cdots .$$
For $Y_0$ we have $d(xy)=0$. Thus, $x dy=-y dx $, which implies $x dx \wedge_R dy=0$ or $y dx \wedge_R dy=0$. But we cannot conclude that $dx \wedge_R dy=0$ on $Y_0$. 

The space of algebraic differential (K$\ddot{a}$hler differential) $1$-forms for $R_0$ is 
$$\Omega^1_{R_0}=\frac{\langle \{dx,dy\} \rangle}{\langle y dx+x dy \rangle}$$
as a $R_0$-module.  Note that for all $n,m \geq 2$, $d^1(x^n dy)=0=d^1(y^m dx)$ on $Y_0$ and $x^n dy, y^m dx$ both are not in $Im$ $d$, but both these elements vanish on $Y_0$. Also, for all $n,m \geq 0$, $x^n dx $ and $y^m dy$ are in $Im$ $d$, provides zero cohomology class.
The K$\ddot{a}$hler differential 1-form $x dy$ (or $y dx$) on $Y_0$ is not in $Im~ d$, but since it is also not in $Ker$ $d^1$ (since $d^1(x dy)= dx \wedge_R dy \neq 0$), does not provides a cohomology class.

Thus, we obtain $H^0_{dR}(Y_0)=Ker$ $d$ $\cong \mathbb{C}$ and $H^1_{dR}(Y_0) \cong \{0\}$.
Next, we compute its second algebraic de Rham cohomology class, for that first consider the $R$-module of algebraic (K$\ddot{a}$hler) differential $2$-forms 
$$\Omega^2_{R_0}=\frac{\langle dx \wedge dy \rangle}{\langle \{x dx\wedge dy, y dx\wedge dy\}\rangle}.$$ Since the only possible choice for a cohomology class in $H^2_{dR}(Y_0)$ is the algebraic $2$-form $dx \wedge dy=d^1(x dy)$ is cohomologous with  zero element,
thus $H^2_{dR}(Y_0)=\{0\}$.  

Note that $Y_0$ is of same homotopy type to a point, thus singular cohomologies of $Y_0$ vanishes in all higher dimensions ($\geq 1$). Fortunately, for this singular variety the algebraic de Rham cohomology and the singular cohomology are same (see \cite{WG}).

This approach does not provide a topological invariant for all singular  algebraic varieties, to resolve this problem  a different approach is considered in \cite{AH}. 

We consider the associated cohomology in complex geometry context (see \cite{AA}).
\begin{Exm}\label{log de Rham cohomology 2} 
 \normalfont{Consider the sheaf of holomorphic functions $\mathcal{O}_X$  on $X:=\mathbb{C}^2$ and the ideal sheaf $\mathcal{I}=\langle xy \rangle$ of $\mathcal{O}_X$. The analytic subspace of $(X, \mathcal{O}_X)$ associated with the ideal sheaf $\mathcal{I}$ is $(Y_0:=V(\mathcal{I}),~ \mathcal{O}_Y:= \mathcal{O}_X/{\mathcal{I}})$.
	Here, the sheaf of logarithmic derivation  and the sheaf of logarithmic $1$-forms are
    $$\mathcal{T}_X(-log~ Y_0)=\langle \{x \partial_x, y \partial_y\} \rangle ~~\text{and}~~\Omega^1_X(log~ Y_0)=\langle \{\frac{dx}{x},\frac{dy}{y}\} \rangle,$$ respectively. Both are locally free $\mathcal{O}_X$-modules of rank $2$.
	The sheaf of logarithmic $2$-forms for the divisor $Y_0$ is $\Omega^2_X(log~ Y_0)=\langle \{\frac{dx}{x}\wedge \frac{dy}{y}\} \rangle$, a locally free $\mathcal{O}_X$-module of rank $1$.
	
	Note that, the sheaf $\mathcal{T}_X(-log~ Y_0)$ has a canonical Lie algebroid structure, and $Y_0$ is a free divisor \cite{BP,AA}. Since,  $log$ $x$ and $log$ $y$ both are not in $\mathcal{O}_X$, thus $\frac{dx}{x}$ and $\frac{dy}{y}$ are not in $\mathscr{I}m$ $d^0$ but both are in $\mathscr{K}er$ $d^1$. Similarly, $\frac{dx}{x}\wedge \frac{dy}{y}$ is in $\mathscr{K}er$ $d^2$ but not in $\mathscr{I}m$ $d^1$. Other choices for cocycles are appears as coboundaries.

	Thus, the hypercohomology of the logarithmic de Rham complex
	\begin{align}\label{log de Rham complex}
		\Omega^{\bullet}_X(log~ Y):	\mathcal{O}_X \overset{d^0}{\longrightarrow}\Omega^1_X(log~ Y_0)\overset{d^1}{\longrightarrow}\wedge^2_{\mathcal{O}_X}\Omega^1_X(log~ Y_0) \overset{d^2}{\longrightarrow}\cdots
	\end{align}
	is $\mathbb{H}^0(X,\Omega^{\bullet}_X(log~ Y_0))=\mathbb{C}$, $\mathbb{H}^1(X,\Omega^{\bullet}_X(log~ Y_0))=\mathbb{C}^2$ and $\mathbb{H}^2(X,\Omega^{\bullet}_X(log~ Y_0))=\mathbb{C}$, all higher cohomologies are zero.
	
	Now, consider the holomorphic Lie algebroid cohomology  \cite{BRT}
	\begin{align}\label{holomorphic Lie algebroid cohomology}
		H^n_{holo}(\mathcal{T}_X(-log~ Y_0)) :=\mathbb{H}^n(X,\Omega^{\bullet}_X(log~ Y_0)) 
		\cong H^n_{dR}(X \setminus Y_0, \mathbb{C})
	\end{align}
	
	(using the standard Lemma of Atiyah-Hodge \cite{MS}), similar results appear in \cite{PT,CMD}. 
	
	Therefore,  the standard Lie algebroid cohomology of $\mathcal{T}_X(-log~ Y_0)$ followed by the isomorphism (\ref{holomorphic Lie algebroid cohomology})  provides the analytic (algebraic) de Rham cohomology of the torus $\mathbb{T}^2$ in $\mathbb{R}^4 \cong \mathbb{C}^2$ (since $X \setminus Y_0= \{(x,y)\in \mathbb{C}^2 \mid x\neq \bar{0}$ and $y\neq \bar{0}\}=(\mathbb{C}\setminus \{\bar{0}\})\times (\mathbb{C}\setminus \{\bar{0}\})$ is a complex manifold and homotopic to $S^1\times S^1= \mathbb{T}^2)$. This agrees with the logarithmic de Rham cohomology for $\Omega^{\bullet}_X(log~ Y_0)$ (described in (\ref{log de Rham complex})) as computed.}
\end{Exm}
\begin{Rem} In general for smooth algebraic variety and complex manifold, to compute its cohomology we need to consider sheaf of  de Rham complex because the space of global sections 
	not necessarily captures the whole information. In these cases we can use \v{C}ech cohomology \cite{VID,SR} by considering an affine open cover or Stein open cover $($known as good open cover$)$ \cite{BRT,MS} accordingly to compute this.
\end{Rem}

\section{Hochschild Cohomology for a Lie algebroid over $a$-spaces} \label{Sec 5}
In this section, we define  Hochschild cohomology of an $\mathcal{O}_X/{\mathbb{K}_X}$-bialgebra (see Section \ref{general bialgebra}), using (Hopf-)Hochschild cohomology of a left bialgebroid \cite{KP} as local descriptions. In particular, we describe the cases associated with  universal enveloping algebroid and jet algebroid of a locally free Lie algebroid over an $a$-space. The cohomology groups can be computed using suitable standard complexes. 
We prove a version of Hochschild-Kostant-Rosenberg (HKR) theorem and dual HKR theorem. It is done by following the local counterpart as can be found in \cite{KP}, where the authors  proved (dual) HKR theorem for (finitely generated) projective Lie-Rinehart algebras. In particular, we present the HKR theorem for the sheaf of logarithmic differential operators and for the sheaf of noncommuative differential operators.

\subsection{Hochschild hypercohomology of an $\mathcal{O}_X/{\mathbb{K}_X}$-bialgebra } \label{Hochschild cohomology of an sheaf of bialgebras}
Let $(\mathcal{A}, \Delta, \epsilon)$ be an $\mathcal{O}_X/{\mathbb{K}_X}$-bialgebra. Consider the presheaf of cochain complexes of $\mathcal{O}_X$-modules
$$U \mapsto C^\bullet(\mathcal{A}(U)),$$ where $C^\bullet(\mathcal{A}(U))=((\mathcal{A}(U))^{\otimes_{\mathcal{O}_X(U)}^\bullet}, b_U)$ is the Hochschild cochain complex of the $\mathcal{O}_X(U)/{\mathbb{K}}$-bialgebra $\mathcal{A}(U)$  
with the differential $b_U: (\mathcal{A}(U))^{\otimes_{\mathcal{O}_X(U)}^\bullet} \rightarrow  (\mathcal{A}(U))^{\otimes_{\mathcal{O}_X(U)}^{\bullet +1}}$ defined for any  $a_1, \dots, a_n \in \mathcal{A}(U)$ as follows:

$$b_U(a_1 \otimes \cdots \otimes a_n)= 1 \otimes a_1 \otimes \cdots \otimes a_n + \sum^n_{i=1} (-1)^i a_1 \otimes \cdots \otimes \Delta (a_i) \otimes \cdots \otimes a_n + (-1)^{n+1}  a_1 \otimes \cdots \otimes a_n \otimes 1~.$$

The sheafification of the presheaf provides a cochain complexes of $\mathcal{O}_X$-modules, we call it as Hochschild cochain complex for $\mathcal{A}$ and denote it by $\mathscr{C}^\bullet(\mathcal{A})$. Its associated hypercohomology is called Hochschild hypercohomology of (the  $\mathcal{O}_X/{\mathbb{K}_X}$-bialgebra) $\mathcal{A}$, and denoted by $\mathbb{H}H^\bullet(\mathcal{A})$. To compute it, it is often useful to consider an appropriate resolution. The canonical choice is given as follows.

Consider the presheaf of cobar complex of left $\mathcal{A}$-comodules
$$U \mapsto  Cob^\bullet({\mathcal{A}(U)})$$ 
where $Cob^\bullet({\mathcal{A}(U)}):=((\mathcal{A}(U))^{\otimes_{\mathcal{O}_X(U)}^{\bullet+1}}, b'_U)$ is the cobar complex for $\mathcal{O}_X(U)$ in the category of left $\mathcal{A}(U)$-comodules
with the differential $b'_U:  (\mathcal{A}(U))^{\otimes_{\mathcal{O}_X(U)}^{\bullet +1}} \rightarrow  (\mathcal{A}(U))^{\otimes_{\mathcal{O}_X(U)}^{\bullet +2}}$ defined for any $a_0, \dots, a_n \in \mathcal{A}(U)$ as follows:
\begin{equation*}\label{cobar}
	b'_U(a_0 \otimes a_1 \otimes \cdots  \otimes a_n)= \sum^n_{i=0}(-1)^i a_0 \otimes \cdots \otimes \Delta(a_i) \otimes \cdots \otimes a_n + (-1)^{n+1} a_0 \otimes \cdots \otimes a_n \otimes 1.
\end{equation*}

The sheafification of the presheaf $U \mapsto Cob^\bullet(\mathcal{A}(U))$, denotes as $\mathscr{C}ob^\bullet(\mathcal{A})$. It provides a resolution of $\mathcal{O}_X$ by cofree  left $\mathcal{A}$-comodules  (a quasi-isomorphism) given by the left $\mathcal{O}_X$-module structure on $\mathcal{A}$ as
$$\mathcal{O}_X \rightarrow \mathscr{C}ob^\bullet(\mathcal{A}),$$ call it cobar resolution of $\mathcal{A}$.
Then by applying cotensor product functor $\mathcal{O}_X \Box_{\mathcal{A}} -$ (sheafifying the cotensor functor described in \cite{KP}) on it and considering hypercohomology (i.e. applying hypercohomology $\mathbb{H}^\bullet(X, -)$ functor on the induced cochain complex $\mathcal{O}_X \Box_{\mathcal{A}}\mathscr{C}ob^\bullet(\mathcal{A})$ of $\mathbb{K}_X$-vector spaces), we get the following isomorphism (sheafifying the local descriptions of relationship  among Hochschild cochain complex and cobar resolution \cite{KP} and considering their associated hypercohomologies)
\begin{align}\label{Hochschild cochain and cobar}
	\mathbb{H}H^\bullet(\mathcal{A}) \cong \mathbb{H}^\bullet (X,\mathcal{O}_X \Box_{\mathcal{A}}\mathscr{C}ob^\bullet(\mathcal{A}))=Cotor^\bullet_{\mathcal{A}}(\mathcal{O}_X,\mathcal{O}_X).
\end{align}

\begin{Rem}
	The standard notion of the Hochschild chain complex is defined for an associative 
	$\mathbb{K}$-algebra using its multiplication, and its dual complex is then considered as the Hochschild cochain complex. In our case, the algebraic structures over $\mathcal{O}_X$ is important to consider. For an $\mathcal{O}_X/{\mathbb{K}_X}$-bialgebra, we do not have an associative algebra structure over $\mathcal{O}_X$, but rather a co-associative algebra structure over $\mathcal{O}_X$. Therefore, for an $\mathcal{O}_X/{\mathbb{K}_X}$-bialgebra the co-associative algebra structure over $\mathcal{O}_X$ defines the Hochschild cohomology using co-multiplication of the $\mathcal{O}_X/{\mathbb{K}_X}$-bialgebra.
\end{Rem}
Now, we describe Hochschild (hyper)cohomology for the $\mathcal{O}_X/{\mathbb{K}_X}$-bialgebras $\mathscr{U}(\mathcal{O}_X,\mathcal{L})$  and $\mathscr{J}(\mathcal{O}_X,\mathcal{L})$ of a (locally free) Lie algebroid $\mathcal{L}$ over an $a$-space $(X,\mathcal{O}_X)$ (see Section \ref{Sec 2} for these notions). 
\subsection{Hochschild hypercohomology of universal enveloping algebroid} \label{Hochschild coho of universal enveloping alg}
First we recall some of the essential  ideas related for our next topic of discussion.

The sheaf of $\mathcal{L}$-poly-vector fields on $X$ is defined as $\mathcal{T}^{poly}_{\mathcal{L}}(\mathcal{O}_X):= \bigoplus_i \wedge^i \mathcal{L}$. This generalizes the space of multisections 
of the tangent bundle (i.e. multi-vector fields)
descried for different $a$-spaces occurs in geometries. 
It has a canonical sheaf of Gerstenhaber algebra structure on $X$.
For the case of a smooth manifold $X$, Kontsevich introduced the (sheaf of) poly-differential operators on $X$. It induces a subcomplex of the Hochschild cochain complex for $\mathcal{O}_X$.
Its analogue in the context of Lie algebroid  is the sheaf of $\mathcal{L}$-poly differential operators $\mathcal{D}^{poly}_{\mathcal{L}}(\mathcal{O}_X):= \bigoplus_i  \mathscr{U}(\mathcal{O}_X, \mathcal{L})^{\otimes_{\mathcal{O}_X}^i}$.  It has a canonical sheaf of Gerstenhaber algebra structure on $X$. In particular, when $\mathcal{L}=\mathcal{T}_X$ we denote these sheaves by the standard notions $\mathcal{T}^{poly}_X$ and $\mathcal{D}^{poly}_X$ respectively.
The so-called HKR map is a quasi-isomorphism between the complexes associated with $\mathcal{T}^{poly}_{X}$ and $\mathcal{D}^{poly}_{X}$ (of differential graded $\mathbb{K}_X$-vector spaces). See \cite{CV} for more generalities.

\label{note} Using ideas from proof of HKR theorem (and canonical PBW coalgebra isomorphism) for Lie-Rinehart algebras \cite{KP} we get the following results. To state these results  first we need to recall some notations associated with a $(\mathbb{K},R)$-Lie-Rinehart algebra $L$ and an $R/{\mathbb{K}}$-bialgebra $A$:

The symmetric algebra of $L$ over $R$ (using underlying $R$-module structure) is denoted by $S_R L$ and  $S_R L^*$ is the symmetric algebra for the $R$-module $L^*:=Hom_R(L,R)$, the universal enveloping algebra of the $(\mathbb{K},R)$-Lie-Rinehart algebra $L$ is denoted by $\mathcal{U}(R,L)$ and its dual is the jet algebra $\mathcal{J}(R,L)$. 

The cobar resolution of an $R/{\mathbb{K}}$-bialgebra $A$ (using underlying $R$-coalgebra structure) is denoted by $Cob^\bullet(A)$ and the associated cohomology (applying the functor $R$ $\Box_A -$ using cotensor product) is Hochschild cohomology of the $R/{\mathbb{K}}$-bialgebra $A$ denoted by $HH^\bullet (A)$ (in our cases $A$ is either $S_RL$, $\mathcal{U}(R,L)$ or $S_R L^*$, $\mathcal{J}(R,L)$ with canonical $R/\mathbb{K}$-bialgebra structure \cite{MM,CV,DRV,KP}).

We know that the sheaf of symmetric algebras $\mathcal{S}_{\mathcal{O}_X}\mathcal{L}$ can be viewed as universal enveloping algebroid of the Lie algebroid $\mathcal{L}$ if the $\mathcal{O}_X$-module $\mathcal{L}$ is equipped with zero bracket and zero anchor (on each space of sections). Recall that we have PBW sheaf homomorphism of  $\mathcal{O}_X$-algebras given as
\begin{center}
	$\theta:\mathcal{S}_{\mathcal{O}_X}\mathcal{L} \rightarrow \mathscr{U}(\mathcal{O}_{X}, \mathcal{L})$\\
	$f \otimes D_1 \otimes \cdots \otimes D_k \mapsto  \frac{1}{k!} f \underset
	{\sigma \in S_k}{\sum}  \bar{D}_{\sigma(1)} \cdots \bar{D}_{\sigma(k)}$
\end{center}
where $f$ is a section of $\mathcal{O}_X$, $D_i$ is a section of $\mathcal{L}$ and $\bar{D}_{i}$ is its associated image $(i=1, 2, \dots,k)$ in $\mathscr{U}(\mathcal{O}_{X}, \mathcal{L})$. Sheafification of the PBW map for Lie-Rinehart algebras (see \cite{GR, JH})  gives the generalized PBW map $\theta$ between the associated sheaves, without considering the graded quotient (see \cite{TS, AA}). 
In addition, if  $\mathcal{L}$ is locally free $\mathcal{O}_X$-module then we get following results.
\begin{Lem}\label{PBW for bialgebra}
	If a Lie algebroid $\mathcal{L}$ is locally free $\mathcal{O}_X$-module, then we have an isomorphism of $\mathcal{O}_X$-coalgebras given by
	the generalized PBW map (or the symmetrization map)$$\theta: \mathcal{S}_{\mathcal{O}_X}\mathcal{L} \rightarrow \mathscr{U}(\mathcal{O}_X,\mathcal{L}).$$
\end{Lem}
\begin{proof}
	First we take
	the PBW map associated with space of sections of $\mathcal{L}$ over each special open subset $U_x$ in $X$ (see Notation \ref{Special open}). It provides an $\mathcal{O}_{X}(U_x)$-coalgebra isomorphism (between the cocommutative $\mathcal{O}_{X}(U_x)$-coalgebras)
	$$ S_{\mathcal{O}_{X}(U_x)}\mathcal{L}(U_x) \cong \mathcal{U}(\mathcal{O}_{X}(U_x), \mathcal{L}(U_x)). $$ It induces stalkwise isomorphisms and by using the sheafification functor over the underlying presheaf homomorphism we get 
	the generalized PBW map $$\mathcal{S}_{\mathcal{O}_X}\mathcal{L} \rightarrow \mathscr{U}(\mathcal{O}_X,\mathcal{L})$$ as an isomorphism of $\mathcal{O}_X$-coalgebras. 
\end{proof}
\begin{Lem}\label{HKR iso for L-R}
	If a Lie algebroid $\mathcal{L}$ is locally free $\mathcal{O}_X$-module, then there is an isomorphism of graded $\mathbb{K}$-vector spaces associated with each special open set $U_x$ (see Notation \ref{Special open}),
	$$HH^\bullet(\mathcal{U}(\mathcal{O}_{X}(U_x), \mathcal{L}(U_x))) \cong \wedge_{\mathcal{O}_{X}(U_x)}^\bullet\mathcal{L}(U_x).$$
\end{Lem}

\begin{proof}
	The $(\mathbb{K},\mathcal{O}_X(U_x))$-Lie-Rinehart algebra $\mathcal{L}(U_x)$ is a free $\mathcal{O}_X(U_x)$-module for each special open set $U_x$ of $X$. As in \cite{KP}, we consider the isomorphism of cochain complexes associated with each special open set $U_x$ given as
	$Cob^\bullet(\mathcal{U}(\mathcal{O}_{X}(U_x), \mathcal{L}(U_x)) \overset{\cong}{\rightarrow} Cob^\bullet(S_{\mathcal{O}_{X}(U_x)}\mathcal{L}(U_x)).$
	
	The isomorphism induces an isomorphism of graded vector spaces $$HH^\bullet(S_{\mathcal{O}_{X}(U_x)}\mathcal{L}(U_x))\cong HH^\bullet(\mathcal{U}(\mathcal{O}_{X}(U_x), \mathcal{L}(U_x)))$$ associated with every special open set $U_x$ for $x\in X$. Now for an open set $U \subset X$ we have the anti-symmetrization (or skew-symmetrization) map \label{Alt}
	
		$$Alt_{U} : \wedge_{\mathcal{O}_{X}(U)}^\bullet\mathcal{L}(U) \rightarrow (S_{\mathcal{O}_{X}(U)}\mathcal{L}(U))^{\otimes^\bullet}$$ 
		 $$D_1 \wedge \cdots \wedge D_n \mapsto \frac{1}{n!}\sum_{\sigma \in S_n} (-1)^\sigma D_{\sigma(1)}\otimes \cdots \otimes D_{\sigma(n)},$$
	as well as the map 
	$$P_{U} : (S_{\mathcal{O}_{X}(U)}\mathcal{L}(U))^{\otimes^\bullet} \rightarrow \wedge_{\mathcal{O}_{X}(U)}^\bullet\mathcal{L}(U)$$ 
		$$\tilde{D_1} \otimes \cdots \otimes \tilde{D_n} \mapsto Pr_{U}(\tilde{D_1})\wedge \cdots \wedge Pr_{U}(\tilde{D_n}),$$
	
	Here $Pr_{U} : S_{\mathcal{O}_{X}(U)}\mathcal{L}(U) \rightarrow S^1_{\mathcal{O}_{X}(U)}\mathcal{L}(U)= \mathcal{L}(U)$ is the projection map on the direct summand $S^1_{\mathcal{O}_X(U)}\mathcal{L}(U)=\mathcal{L}(U)=\wedge^1_{\mathcal{O}_X(U)} \mathcal{L}(U)$, and $D_i \in \mathcal{L}(U)$, $\tilde{D_i} \in S_{\mathcal{O}_{X}(U)}\mathcal{L}(U)$ for $i= 1,\dots ,n$.
	
	These morphisms defines cochain equivalence 
	$Cob^\bullet(S_{\mathcal{O}_X(U_x)}\mathcal{L}(U_x))\overset{\sim}{\rightarrow} (\wedge_{\mathcal{O}_{X}(U_x)}^\bullet\mathcal{L}(U_x),0)$ (i.e. they form a quasi isomorphism) 
	on each special open set $U_x$ in $X$ and using that we get the isomorphism of graded vector spaces
	\begin{center}
		$HH^\bullet(\mathcal{U}(\mathcal{O}_{X}(U_x), \mathcal{L}(U_x))) \cong \wedge_{\mathcal{O}_{X}(U_x)}^\bullet\mathcal{L}(U_x)$.
	\end{center}
\end{proof}
\vspace{-.1 cm}

On each open set $U$ of $X$, we can define the anti-symmetrization map 
$$ \widetilde{Alt}_U: \wedge^\bullet_{\mathcal{O}_{X}(U)}\mathcal{L}(U) \rightarrow (\mathcal{U}(\mathcal{O}_X(U),\mathcal{L}(U)))^{\otimes^\bullet}$$
as the earlier case and we can check it is compatible with restrictions and differentials. Applying sheafification functor, we get homomorphism between associated cochain complex of $\mathcal{O}_X$-modules, we call it anti-symmetrization map (it depends on context). We use these results in the next theorem.

\begin{Lem} \label{hyper & sheaf}
	If $\mathcal{F}$ is a sheaf of $\mathcal{O}_X$-modules, then by considering the exterior algebra as a chain complex $\wedge^\bullet _{\mathcal{O}_X} \mathcal{F}:=(\oplus_i \wedge^i _{\mathcal{O}_X} \mathcal{F},0)$ we get a canonical isomorphism of $\mathbb{K}$-vector spaces
	$$\mathbb{H}^k(X, \wedge^\bullet _{\mathcal{O}_X} \mathcal{F}) \cong \underset{i+j=k}{\oplus} H^j(X, \wedge^i _{\mathcal{O}_X} \mathcal{F})$$ for every $k \geq 0$ (or $k$ is any non negative integer).
\end{Lem}
\begin{proof} Here we apply the notion of hypercohomology (see Section \ref{Hypercohomology}) for the cochain complex of sheaves $\wedge^\bullet _{\mathcal{O}_X} \mathcal{F}$.
	For each sheaf of ${\mathcal{O}_X}$-modules (or $\mathcal{O}_X$-module) $\wedge^i _{\mathcal{O}_X} \mathcal{F}$, we have an injective resolution (given by the flabby Godement resolution) $\mathscr{C}^\bullet(\wedge^i _{\mathcal{O}_X} \mathcal{F})$, provides a quasi-isomorphism
	$$\wedge^i _{\mathcal{O}_X} \mathcal{F} \overset{\sim}{\rightarrow} \mathscr{C}^\bullet(\wedge^i _{\mathcal{O}_X} \mathcal{F}).$$
	Thus, the sheaf cohomology $H^\bullet(X,\wedge^i _{\mathcal{O}_X} \mathcal{F})$ of $\wedge^i _{\mathcal{O}_X} \mathcal{F}$ ($i \in \mathbb{N} \cup \{0\})$ is isomorphic to the cohomology (in usual sense) of the complex of $\mathcal{O}_X(X)$-modules $\mathscr{C}^\bullet(\wedge^i _{\mathcal{O}_X} \mathcal{F})(X)$.

	The corresponding bicomplex for the complex of sheaves $\wedge^\bullet _{\mathcal{O}_X} \mathcal{F}$ is $\mathscr{C}^\bullet(\mathcal{F}):=(\mathscr{C}^j(\wedge^i _{\mathcal{O}_X} \mathcal{F}))_{i,j \geq 0}$. The original complex  can be embedded in the total complex $\mathscr{K}^\bullet:=tot(\mathscr{C}^\bullet(\mathcal{F}))$. Moreover, this embedding is a quasi-isomorphism. The cohomology of the associated complex of global sections $\mathscr{K}^\bullet(X)=tot(\mathscr{C}^\bullet(\mathcal{F}))(X)$ is the hypercohomology of the complex $\wedge^\bullet _{\mathcal{O}_X} \mathcal{F}$. Since here differential of the complex is zero, thus the computation of cohomology of the total complex is become a simpler one. It is directly expressed through the sheaf cohomologies $H^j(X, \wedge^i _{\mathcal{O}_X} \mathcal{F})$ as stated in the Lemma.
	
\end{proof}
\begin{Nota}
For an open set $U$ of $X$, the cotensor product with $\mathcal{O}_X(U)$ in the category of left $S_{\mathcal{O}_X(U)}\mathcal{L}(U)$-comodules is  $\mathcal{O}_X(U)$ $\Box_{S_{\mathcal{O}_X(U)}\mathcal{L}(U)}$ $-$ and the associated cohomology groups is given by the Cotor groups (functor) as $Cotor^\bullet_{S_{\mathcal{O}_X(U)}\mathcal{L}(U)} (\mathcal{O}_X(U),-)$. Instead of $S_{\mathcal{O}_X(U)}\mathcal{L}(U)$, we use $\mathcal{J}(\mathcal{O}_X(U), \mathcal{L}(U))$ when it requires.
\end{Nota}
\begin{Thm} $($Generalized HKR theorem$)$ \label{Generalized HKR thm} Let $\mathcal{L}$ be a locally free  Lie algebroid over $(X, \mathcal{O}_X)$. Then the anti-symmetrization map (known as HKR morphism)
	\begin{center}
		$\wedge^\bullet_{\mathcal{O}_X} \mathcal{L} \rightarrow \mathscr{U}(\mathcal{O}_X, \mathcal{L})^{\otimes^ \bullet}$
	\end{center}
	provides a quasi-isomorphism of the associated cochain complexes of sheaves. It induces an isomorphism of graded vector spaces 
	\begin{center}
		$ \mathbb{H}H^\bullet(\mathscr{U}(\mathcal{O}_X, \mathcal{L})) \cong \mathbb{H}^\bullet(X, \wedge^\bullet_{\mathcal{O}_X}{\mathcal{L}})$
	\end{center}
\end{Thm}

\begin{proof} Assume first $\mathcal{L}$ to be a locally free Lie algebroid of finite rank. Since for each $x \in X$ we have an open set $U_x$ containing $x$ with $\mathcal{L}|_{U_x}$ is a free $\mathcal{O}_X|_{U_x}$-module, thus on each $U_x$ we use results from the proof of the HKR Theorem \cite{KP} for the $(\mathbb{K}, \mathcal{O}_X(U_x))$-Lie-Rinehart algebra $\mathcal{L}(U_x)$.
	Then, for each special open sets $U_x$ (see Notation \ref{Special open}) we have isomorphisms of cochain complexes of left $ S_{\mathcal{O}_X(U_x)}\mathcal{L}(U_x)$-comodules
	\begin{center}
		$ (\mathcal{O}_X (U_x)$ $\Box_{S_{\mathcal{O}_X(U_x)}\mathcal{L}(U_x)}$ $ \tilde{K}^\bullet_{U_x}\mathcal{L},$ $ id_{\mathcal{O}_X(U_x)}$ $\otimes_{\mathcal{O}_X(U_x)}$ $ \partial_{U_x}) \overset{\phi_{U_x}}{\rightarrow} (\wedge^\bullet_{\mathcal{O}_X(U_x)} \mathcal{L}(U_x), 0)$
	\end{center}
	compatible with restrictions, where $\tilde{K}^\bullet_{U_x}\mathcal{L} := S_{\mathcal{O}_X(U_x)}\mathcal{L}(U_x)\otimes_{\mathcal{O}_X(U_x)} \wedge^\bullet_{\mathcal{O}_X(U_x)}\mathcal{L}(U_x)$ is with the canonical $S_{\mathcal{O}_X(U_x)}\mathcal{L}(U_x)$-comodule structure and the differential $\partial_{U_x}: \tilde{K}^\bullet_{U_x} \mathcal{L} \rightarrow \tilde{K}^{\bullet+1}_{U_x} \mathcal{L}$ is given by the dual Koszul resolution of $\mathcal{O}_X(U_x)$ as $S_{\mathcal{O}_X(U_x)}\mathcal{L}(U_x)$-comodules.
	Indeed, the dual Koszul-Rinehart complex for  $\mathcal{L}(U)$ is given as $\tilde{K}^\bullet_{\mathcal{O}_X(U)}\mathcal{L}(U):=(\tilde{K}^\bullet_{U}\mathcal{L} , $ $\partial_U)$, for any open set $U$ in $X$ $($see Remark \ref{dual K-R}$)$.

	Now, an open set $U \subset X$ is expressed by taking a partition as $\cup_{x \in X}(U \cap U_x)$ and thus we have cochain homomorphisms of left $S_{\mathcal{O}_X(U)}\mathcal{L}(U)$-comodules
	\begin{center} \label{Koszul and wedge}
		$ (\mathcal{O}_X (U)$ $\Box_{S_{\mathcal{O}_X(U)}\mathcal{L}(U)}$ $ \tilde{K}^\bullet_{U}\mathcal{L},$ $ id_{\mathcal{O}_X(U)}$ $\otimes_{\mathcal{O}_X(U)} $ $ \partial_{U}) \overset{\phi_{U}}{\rightarrow} (\wedge^\bullet_{\mathcal{O}_X(U)} \mathcal{L}(U), 0)$
	\end{center}
	compatible with restrictions. Then it induces a homomorphism on the associated presheaves of cochain complexes, provides stalkwise isomorphism. Thus, applying sheafification functor we get a canonical isomorphism of cochain complexes of left $\mathcal{S}_{\mathcal{O}_X}\mathcal{L}$-comodules $($see Remark
    \ref{dual K-R}$)$
	\begin{align} \label{Koszul 1}
		\mathcal{O}_X \Box_{\mathcal{S}_{\mathcal{O}_X}\mathcal{L}}~  \widetilde{\mathscr{K}}^\bullet_{\mathcal{O}_X}\mathcal{L}=(\mathcal{O}_X \Box_{\mathcal{S}_{\mathcal{O}_X}\mathcal{L}}~ (\mathcal{S}_{\mathcal{O}_X}\mathcal{L}\otimes \wedge^\bullet_{\mathcal{O}_X}\mathcal{L}), id_{\mathcal{O}_X} \otimes_{\mathcal{O}_X} \partial) \overset{\phi}{\rightarrow} (\wedge^\bullet_{\mathcal{O}_X} \mathcal{L}, 0).
	\end{align}

	To compute Hochschild hypercohomology $\mathbb{H}H^\bullet(\mathcal{S}_{\mathcal{O}_X} \mathcal{L})$ via a derived functor, we consider the standard cobar complex associated with each open set $U$, denoted as $Cob^\bullet(S_{\mathcal{O}_X(U)} \mathcal{L}(U))$ (considering $\mathcal{A}(U)=S_{\mathcal{O}_X(U)} \mathcal{L}(U)$ as described in Section \ref{Hochschild cohomology of an sheaf of bialgebras}).
	
	We have two presheaf of vector spaces, one is presheaf of cobar complexes of the sheaf $\mathcal{S}_{\mathcal{O}_X} \mathcal{L}$ of $\mathcal{O}_X/{\mathbb{K}_X}$-bialgebra, defined as
	$$U \mapsto Cob^\bullet(S_{\mathcal{O}_X(U)} \mathcal{L}(U))$$
	with its restriction morphism induced from the restrictions of the sheaf $\mathcal{S}_{\mathcal{O}_X}\mathcal{L}$, and another one is the presheaf of Koszul-Rinehart complex of $\mathcal{O}_X$ in the category of left $S_{\mathcal{O}_X}\mathcal{L}$-comodules, defined as
	$$U\mapsto \widetilde{K}_{\mathcal{O}_X(U)}^\bullet\mathcal{L}(U).$$ 
	These two presheaves are quasi isomorphic on special open sets $U_x$ for every $x \in X$ (for the $C^\infty$ case these holds for every open set $U$). 
	In stalkwise, it induces quasi-isomorphism and by applying sheafification functor we get quasi-isomorphism (chain homotopy) 
	\begin{align} \label{Cob and Kos}
		\mathscr{C}ob^\bullet(\mathcal{S}_{\mathcal{O}_X} \mathcal{L}) \overset{\sim}{\rightarrow} \widetilde{\mathscr{K}}_{\mathcal{O}_X}^\bullet\mathcal{L}
	\end{align}
	(using the alternating map and projection map (\ref{Alt}), we have quasi isomorphism $$S_{\mathcal{O}_X(U_x)}\mathcal{L}(U_x) \otimes_{\mathcal{O}_X(U_x)} \wedge^\bullet_{\mathcal{O}_X(U_x)} (\mathcal{L}(U_x)) \overset{id \otimes Alt_{U_x}}{\longrightarrow} S(\mathcal{O}_X(U_x),\mathcal{L}(U_x))^{\otimes^{(\bullet+1)}}).$$
	Then on taking cotensor product by $\mathcal{O}_X$ from left (as standard left $\mathcal{S}_{\mathcal{O}_X}\mathcal{L}$-comodule) to the cochain complexes (\ref{Cob and Kos}) of sheaves of left $\mathcal{S}_{\mathcal{O}_X}\mathcal{L}$-comodules we get from (\ref{Koszul 1})
	\begin{center}
		$\mathcal{O}_X \Box_{\mathcal{S}_{\mathcal{O}_X} \mathcal{L}}$ $ \mathscr{C}ob^\bullet(\mathcal{S}_{\mathcal{O}_X} \mathcal{L}) \overset{\sim}{\rightarrow} (\wedge^\bullet_{\mathcal{O}_X}\mathcal{L}, 0)$.
	\end{center}
	Applying hypercohomology functor $\mathbb{H}^\bullet(X,-)$, we get Hochschild hypercohomology of $\mathscr{U}(\mathcal{O}_X,\mathcal{L})$ (using Lemma \ref{PBW for bialgebra}) through derived functor Cotor as
	\begin{center}
		$\mathbb{H}H^\bullet(\mathscr{U}(\mathcal{O}_X,\mathcal{L})) \cong Cotor^\bullet_{\mathcal{S}_{\mathcal{O}_X}\mathcal{L}}(\mathcal{O}_X, \mathcal{O}_X) \overset{\cong}  {\rightarrow} \mathbb{H}^\bullet(X, \wedge^\bullet_{\mathcal{O}_X}\mathcal{L})$.
	\end{center}
	An alternative approach for proof of the theorem for finite rank case is given as follows.
	
	Consider the  presheaves of graded vector spaces
	\begin{center}
		$U \mapsto HH^\bullet(S_{\mathcal{O}_{X}(U)}\mathcal{L}(U))$, 
		\vspace{.2 cm}
		
		\hspace{-.7 cm}		$U \mapsto \wedge_{\mathcal{O}_{X}(U)}^\bullet$ $\mathcal{L}(U)$,
	\end{center}
	which are isomorphic on each special open sets $U_x$'s (see Lemma \ref{HKR iso for L-R}). The associated presheaf homomorphism is induced from the homomorphisms $Alt_U$ and $P_U$ (described in the proof of the Lemma \ref{HKR iso for L-R}).
	
	Thus, by considering sheafifications (of these two presheaf of cohomology spaces which are isomorphic stalkwise) we get an isomorphism in the associated cohomology sheaves (described in Section \ref{Hypercohomology}) as
	$$\mathscr{H}H^\bullet(\mathscr{U}(\mathcal{O}_{X}, \mathcal{L}))\cong \mathscr{H}H^\bullet(\mathcal{S}_{\mathcal{O}_X}\mathcal{L}) \cong \wedge_{\mathcal{O}_X}^\bullet\mathcal{L}.$$
	
	Next on the induced hypercohomology groups (using Section \ref{Hypercohomology}) we have
	\begin{center}
		$\mathbb{H}H^\bullet(\mathscr{U}(\mathcal{O}_{X}, \mathcal{L}))\cong \mathbb{H}H^\bullet(\mathcal{S}_{\mathcal{O}_X}\mathcal{L}) \cong \mathbb{H}^\bullet (X, \wedge_{\mathcal{O}_X}^\bullet\mathcal{L})$.
	\end{center}
	In general case, where $\mathcal{L}$ is locally free $\mathcal{O}_X$-module of infinite rank (quasicoherent sheaf), there exists a filtered ordered set $J$ as well as an inductive system of free $\mathcal{O}_X(U_x)$-modules $\{{\mathcal{L}(U_x)}_j \mid j \in J \}$ such that $$\mathcal{L}(U_x) \cong \varinjlim_{ j \in J} {\mathcal{L}(U_x)}_j,$$ for each special open set $U_x$  (using results from \cite{KP}). Since both $\mathbb{H}H^\bullet$ (which is isomorphic to the derived functor $Cotor^\bullet$) and the functor $\mathcal{S}$ commute with inductive limits over a filtered ordered set, we get the result using the result follows from the locally free $\mathcal{O}_X$-module of finite rank case (local counterpart is described in \cite{KP}).
\end{proof}
\begin{Rem}
	The HKR morphism is not a canonical ring isomorphism, but composing it together with the Todd genus provides a ring isomorphism (moreover, forms a canonical Gerstenhaber algebra isomorphism). It is required in the study of formality (or quantization) for Lie algebroids over a ringed site \cite{CV}.
\end{Rem}
\begin{Cor} \label{simplify HKR} By applying the Lemma \ref{hyper & sheaf}, the above HKR theorem reduces to
	$$ \mathbb{H}H^\bullet(\mathscr{U}(\mathcal{O}_X, \mathcal{L}))\cong \underset{i,j}{\oplus} H^j(X, \wedge^i _{\mathcal{O}_X} \mathcal{L}).$$
\end{Cor}
\begin{Cor}
	In the special case of $\mathcal{L} = \mathcal{T}_X $ when $X$ is non-singular or smooth  $($see Remark \ref{special $a$-spaces}$)$,
	we get isomorphism of graded vector spaces 
	\begin{center}
		$ \mathbb{H}H^\bullet(\mathcal{D}_X)= \mathbb{H}^\bullet(X, \mathcal{D}^{poly}_X) \cong \mathbb{H}^\bullet(X, \mathcal{T}^{poly}_X)$.
	\end{center}
\end{Cor}
\begin{Cor}
	Using Example \ref{sheaf of log diff operators} for a free divisor $Y$ in $X$, we get isomorphism of graded vector spaces 
	$$\mathbb{H}H^\bullet(\mathcal{D}_X(-log Y)) \cong \mathbb{H}^\bullet(X, \wedge^\bullet_{\mathcal{O}_X} \mathcal{T}_X(-log Y)).$$
\end{Cor}
\begin{Cor} $($Noncommutative analogue of the HKR Theorem$)$
	The free Lie algebroid $\mathcal{P}_X$ can be identified with sheaf of regular noncommutative vector fields on $X$ and  we have seen already that the sections of its universal enveloping algebra $\mathbb{D}_X$ are noncommutative differential operators on $X$ $($see Example \ref{free Lie algebroid}$)$. Thus, in this case from the HKR theorem we get a canonical isomorphism
	$$\mathbb{H}H^\bullet(\mathbb{D}_X) \cong \mathbb{H}^\bullet(X,\wedge^\bullet_{\mathcal{O}_X}\mathcal{P}_X)$$
\end{Cor}	

\begin{Cor}
	For  a cocomplete graded free $\mathcal{O}_X/{\mathbb{K}_X}$-bialgebra $\mathcal{A}$ $($of finite or infinite type$)$, the isomorphism $\mathcal{A}\cong \mathscr{U}(\mathcal{O}_X,\mathscr{P}(\mathcal{A}))$ holds \cite{AA} and thus we get the graded vector space isomorphism
	$$\mathbb{H}H^\bullet(\mathcal{A}) \cong \mathbb{H}^\bullet(X, \wedge^\bullet_{{\mathcal{O}_X}} \mathscr{P}(\mathcal{A})).$$
\end{Cor}
\subsection{Hochschild hypercohomology of jet algebroid}\label{Hochschild coho of jet algebroid} Here, we present a dual version of the HKR theorem by considering Hochschild cohomology of the jet algebroid $\mathscr{J}(\mathcal{O}_X,\mathcal{L})$ for a Lie algebroid $\mathcal{L}$ over $(X,\mathcal{O}_X)$ (see Section \ref{jet alg}).
\begin{Rem}\label{Grothendick conn}
	There is a canonical left $\mathscr{U}(\mathcal{O}_X,\mathcal{L})$-module structure on $\mathscr{J}(\mathcal{O}_X,\mathcal{L})$ constructed as follows $($see \cite{CV,KP,CRV} for the local descriptions$)$: 
	
	A canonical flat $\mathcal{L}$-connection on $\mathscr{J}(\mathcal{O}_X,\mathcal{L})$, called Grothendieck connection is given by
	\begin{align} \label{Grothendick connection}
		\nabla_{{D}} (\phi)(D') := \tilde{D}(\phi(D')) - \phi(\bar{D}~D'),
	\end{align}
	for all sections ${D}\in \mathcal{L}$, $\phi \in \mathscr{J}(\mathcal{O}_X,\mathcal{L})$ and $D' \in \mathscr{U}(\mathcal{O}_X,\mathcal{L})$, where $\tilde{D}:=\mathfrak{a}(D)$ and $\bar{D}:=\iota_{\mathcal{L}}(D)$. 
	Thus, we get an $\mathcal{L}$-module structure on $\mathscr{J}(\mathcal{O}_X,\mathcal{L})$,  can be extended to 
	an induced $\mathscr{U}(\mathcal{O}_X,\mathcal{L})$-module structure on $\mathscr{J}(\mathcal{O}_X,\mathcal{L})$.
\end{Rem}

\begin{Thm} $($Generalized dual HKR theorem$)$ \label{Generalized dual HKR thm} Let $\mathcal{L}$ be a locally free Lie algebroid over $(X, \mathcal{O}_X)$ which is of finite rank. Then the anti-symmetrization map $($dual HKR morphism$)$
	\begin{center}
		$\wedge^\bullet_{\mathcal{O}_X} \mathcal{L^*} \rightarrow \mathscr{J}(\mathcal{O}_X, \mathcal{L})^{\otimes^ \bullet}$
	\end{center}
	provides a quasi-isomorphism between the associated cochain complexes of sheaves. It induces canonical isomorphism of graded vector spaces 
	\begin{center}
		$ \mathbb{H}H^\bullet(\mathscr{J}(\mathcal{O}_X,\mathcal{L})) \cong \mathbb{H}^\bullet(X, \Omega^\bullet_{\mathcal{L}}).$
	\end{center}
	
\end{Thm}

\begin{proof} 
	Consider a locally free Lie algebroid $\mathcal{L}$ over $(X, \mathcal{O}_X)$ of finite rank, say $r$ for some $r \in \mathbb{N}$, and the dual of $\mathcal{L}$, which is denoted by $\mathcal{L^*}$, is a locally free $\mathcal{O}_X$-module of the same rank. Thus, $\mathcal{L^*}(U_x)$ has a basis $\{w_1, \dots, w_r\}$ (say) for each special open set $U_x$ around $x \in X$. Since $\mathscr{J}(\mathcal{O}_X,\mathcal{L})$ is a commutative $\mathcal{O}_X$-algebra, thus we obtain a canonical $\mathcal{O}_X(U_x)$-algebra isomorphism
    $$\mathcal{J}(\mathcal{O}_X(U_x),\mathcal{L}(U_x))
	\cong \mathcal{O}_X (U_x)[[w_1, \dots, w_r]].$$  Hence, by applying sheafification functor we get the isomorphism of $\mathcal{O}_X$-algebras as $\mathscr{J}(\mathcal{O}_X,\mathcal{L}) \cong \widehat{\mathcal{S}_{\mathcal{O}_X}\mathcal{L^*}}$ (this sheaf of symmetric algebras is formally completed with respect to the degree) \cite{DRV,AP}. 
    
    Now, consider the dual Koszul-Rinehart resolution of $\mathcal{O}_X$ in the category of left $\mathscr{J}(\mathcal{O}_X,\mathcal{L})$-comodules (by using local counterpart from \cite{KP}) is
	\begin{center}
		$ \mathcal{O}_X   \hookrightarrow (\mathscr{J}(\mathcal{O}_X,\mathcal{L}) \otimes_{\mathcal{O}_X} \wedge^\bullet_{\mathcal{O}_X}\mathcal{L^*}, \bar{\nabla}) = : \bar{\mathscr{K}}^\bullet_{\mathcal{O}_X}\mathcal{L}^*,$
	\end{center}
	where $\bar{\nabla}$ is the continuation of the Grothendieck connection $\nabla$ (Remark \ref{Grothendick conn}), where $\nabla$ is viewed as (defined in (\ref{L-connection special}))   $\nabla: \mathscr{J}(\mathcal{O}_X,\mathcal{L}) \rightarrow \mathscr{J}(\mathcal{O}_X,\mathcal{L}) \otimes_{\mathcal{O}_X} \mathcal{L^*}$. 
	  It is basically the sheafification of the presheaves of cochain complex
	$$U \mapsto \{\mathcal{O}_X(U) \hookrightarrow \bar{K}^\bullet_{\mathcal{O}_X(U)} \mathcal{L}(U)^* \},$$
	where $\bar{K}^\bullet_{\mathcal{O}_X(U)} \mathcal{L}(U)^*= (\mathcal{J}(\mathcal{O}_X(U),\mathcal{L}(U)) \otimes_{\mathcal{O}_X(U)} \wedge^\bullet_{\mathcal{O}_X(U)}\mathcal{L}(U)^*, \bar{\nabla}^U)$ and $\bar{\nabla}: U \mapsto \bar{\nabla}^U$ is the sheafification. The unit of $\mathscr{J}(\mathcal{O}_X,\mathcal{L})$ as an $\mathcal{O}_X$-algebra provides the morphism from $\mathcal{O}_X$ to $\bar{\mathscr{K}}^\bullet_{\mathcal{O}_X}\mathcal{L}^*$.
	
	Applying cotensor product by $\mathcal{O}_X$ (with canonical $\mathscr{J}(\mathcal{O}_X,\mathcal{L})$-comodule structure, which induces from the $\mathscr{U}(\mathcal{O}_X,\mathcal{L})$-module structure on $\mathcal{O}_X$) to the dual Koszul-Rinehart complex $\bar{\mathscr{K}}^\bullet_{\mathcal{O}_X}\mathcal{L}^*$, we get the canonical isomorphism (sheafifying local descriptions from \cite{KP}) of cochain complexes of sheaves of  graded vector spaces as
	\begin{align} \label{Iso}
		\mathcal{O}_X \Box_{\mathscr{J}(\mathcal{O}_X,\mathcal{L})} \bar{\mathscr{K}}^\bullet_{\mathcal{O}_X}\mathcal{L}^* \cong (\wedge^\bullet_{\mathcal{O}_X} \mathcal{L^*}, d).
	\end{align}
	To describe the isomorphism (\ref{Iso}) in an explicit way, we use the following steps.
	$$ \mathcal{O}_X (U)\Box_{\mathcal{J}(\mathcal{O}_X(U),\mathcal{L}(U))} (\mathcal{J}(\mathcal{O}_X(U),\mathcal{L}(U)) \otimes_{\mathcal{O}_X(U)}\wedge^\bullet_{\mathcal{O}_X(U)} \mathcal{L}(U)^*) \cong\wedge^\bullet_{\mathcal{O}_X(U)} \mathcal{L}(U)^*.$$
	Since the unit $1_U \in \mathcal{J}(\mathcal{O}_X(U),\mathcal{L}(U))$ is given by the  counit $\epsilon_U$ of $\mathcal{U}(\mathcal{O}_X(U),\mathcal{L}(U))$, the induced differential is exactly the Lie-Rinehart coboundary $d_{U}$ (\ref{Chevalley-Eilenberg-de Rham complex}), on an open set $U \subset X$.
	These isomorphisms  are compatible with restrictions.
	Now consider the canonical presheaves from the above complexes associated with each open sets, which are isomorphic on each special open sets (see Notation \ref{Special open}) and compatible with the natural restrictions. Then considering there sheafifications and using the above isomorphisms, we get isomorphism between the complexes of sheaves, described in (\ref{Iso}).
	
	Since $\mathcal{L}$ is locally free $\mathcal{O}_X$-module of finite rank, there exists open sets $U_x$ around each point $x\in X$, where $\mathcal{L}|_{U_x}$ is free $\mathcal{O}_X|_{U_x}$-module of finite rank. Thus, for each $U_x$ we can use results from the proof of the dual HKR Theorem \cite{KP} for the $(\mathbb{K}, \mathcal{O}_X(U_x))$-Lie-Rinehart algebra $\mathcal{L}(U_x)$. Here, we get
	$$\wedge^\bullet_{\mathcal{O}_X(U_x)} \mathcal{L}(U_x)^* \cong (\wedge^\bullet_{\mathcal{O}_X(U_x)} \mathcal{L}(U_x))^*.$$
	Hence, $\mathscr{H}om_{\mathcal{O}_X}(\wedge^\bullet_{\mathcal{O}_X} \mathcal{L}, \mathcal{O}_X)\cong \wedge^\bullet_{\mathcal{O}_X} \mathcal{L^*}$ and the Chevalley-Eilenberg-de Rham complex of $\mathcal{L}$ is $ \Omega^\bullet_{\mathcal{L}} \cong (\wedge^\bullet_{\mathcal{O}_X} \mathcal{L^*}, d)$. Therefore, we get (using isomorphism (\ref{Iso}))
	\begin{align}\label{Iso 2}
		\mathcal{O}_X \Box_{\mathscr{J}(\mathcal{O}_X,\mathcal{L})} \bar{\mathscr{K}}^\bullet_{\mathcal{O}_X}\mathcal{L}^* \cong \Omega^\bullet_{\mathcal{L}}.
	\end{align}
	
	To express Hochschild hypercohomology of $\mathscr{J}(\mathcal{O}_X, \mathcal{L})$ as derived functor, we need to consider standard cobar resolution of $\mathcal{O}_X$ as $\mathscr{J}(\mathcal{O}_X, \mathcal{L})$-comodules and cotensor it with $\mathcal{O}_X$ (putting $\mathcal{A}= \mathscr{J}(\mathcal{O}_X,\mathcal{L})$ in the relation appears in (\ref{Hochschild cochain and cobar})). 
	Thus, we obtain the graded vector space isomorphism
	\begin{align} \label{Cochain,Cobar}
		\mathbb{H}H^\bullet(\mathscr{J}(\mathcal{O}_X,\mathcal{L})) \cong Cotor^\bullet_{\mathscr{J}(\mathcal{O}_X, \mathcal{L})}(\mathcal{O}_X, \mathcal{O}_X).
	\end{align}
	
	Since for each open set $U_x$ of $X$, the $(\mathbb{K}, \mathcal{O}_X(U_x))$-Lie-Rinehart algebras $\mathcal{L}(U_x)$ is finitely generated projective (in fact free), thus we get quasi-isomorphisms  
	$$ \mathcal{J}(\mathcal{O}_X(U_x),\mathcal{L}(U_x))^{\otimes (\bullet +1)} \overset{id_{U_x} \otimes P_{U_x}}{\longrightarrow} \mathcal{J}(\mathcal{O}_X(U_x),\mathcal{L}(U_x)) \otimes_{\mathcal{O}_X(U_x)} \wedge^\bullet_{\mathcal{O}_X(U_x)} \mathcal{L}(U_x)^*,$$
	
	$$\mathcal{J}(\mathcal{O}_X(U_x),\mathcal{L}(U_x)) \otimes_{\mathcal{O}_X(U_x)} \wedge^\bullet_{\mathcal{O}_X(U_x)} \mathcal{L}(U_x)^* \overset{id_{U_x} \otimes Alt_{U_x}}{\longrightarrow} (\mathcal{J}(\mathcal{O}_X(U_x),\mathcal{L}(U_x)))^{\otimes (\bullet +1)},$$
	for each $U_x$ associated with $x \in X$, where $P_{U_x} $ and $Alt_{U_x}$ is given by the canonical projections and the anti-symmetrization map, respectively, given as
	$$pr_1 : \mathcal{J}(\mathcal{O}_X(U_x),\mathcal{L}(U_x)) \cong \widehat{S_{\mathcal{O}_X (U_x)}\mathcal{L}(U_x)^*} \rightarrow \mathcal{L}(U_x)^*$$  
	$$Alt_{U_x} : \wedge^\bullet_{\mathcal{O}_X(U_x)} \mathcal{L}(U_x)^* \rightarrow \mathcal{J}(\mathcal{O}_X(U_x),\mathcal{L}(U_x))^{\otimes^ \bullet} $$ for every special open set $U_x$ of $ X$ (for smooth manifold we get these quasi-isomorphisms on each open sets).
	Sheafification of the associated presheaves provides the quasi-isomorphism between the cobar complex $\mathscr{C}ob^\bullet(\mathscr{J}(\mathcal{O}_X, \mathcal{L}))$ and the dual Koszul-Rinehart complex $\bar{\mathscr{K}}^\bullet_{\mathcal{O}_X}\mathcal{L}^*$.
	
	Applying hypercohomology functor $\mathbb{H}^\bullet(X,-)$ we get (from the isomorphism (\ref{Iso 2}))
	\begin{align} \label{Cotor, de Rham}
		Cotor^\bullet_{\mathscr{J}(\mathcal{O}_X, \mathcal{L})}(\mathcal{O}_X, \mathcal{O}_X) \cong \mathbb{H}^\bullet(X, \Omega^\bullet_{\mathcal{L}}).
	\end{align}
	
	Thus, using the isomorphisms (\ref{Cochain,Cobar}) and (\ref{Cotor, de Rham}), the Hochschild hypercohomology groups of jet algebroid $\mathscr{J}(\mathcal{O}_X, \mathcal{L})$ is expressed by the Chevally-Eilenberg-de Rham hypercohomology of $\mathcal{L}$.
\end{proof}
\begin{Rem}
	In \cite{DRV,CRV}, a dual version of the (twisted) HKR theorem for Lie algebroids \cite{CV}, utilizing different techniques to handle rich algebraic structures, is considered to study the notion of precalculus up to homotopy in the context of Lie algebroids over ringed sites.
\end{Rem}
\begin{Cor} \label{Hyper coho of jets and de Rham coho}
	In particular, when $\mathcal{L} = \mathcal{T}_X $ over a non-singular $a$-space (smooth manifold, complex manifold, smooth algebraic variety or smooth scheme over the field $\mathbb{C})$ $X$,
	we get isomorphism of graded vector spaces 
	\begin{center}
		$ \mathbb{H}H^\bullet(\mathscr{J}(\mathcal{O}_X,\mathcal{T}_X)) \cong \mathbb{H}^\bullet(X, \Omega^\bullet_X) \cong H^\bullet(X, \mathbb{K})$ or $H^\bullet(X^{an}, \mathbb{K})$ 
	\end{center}
	by applying de Rham theorems in different settings $($smooth, analytic, algebraic$)$ \cite{VD,MS},
	where $ \mathbb{K = R}$ or $ \mathbb{C}$ accordingly and $X^{an}$ is the analytification of the algebraic variety (scheme) $X$.
\end{Cor}
\begin{Cor}
	Applying the Theorem \ref{Lie algebroid Cohomology as derived functor} for $\mathcal{L} = \mathcal{T}_X $ over some non-singular $a$-space $X$ (see Remark \ref{special $a$-spaces}) and using the Corollary \ref{Hyper coho of jets and de Rham coho} we get the isomorphism
	$$ \mathbb{H}H^\bullet(\mathscr{J}(\mathcal{O}_X,\mathcal{T}_X)) \cong Ext^\bullet_{\mathcal{D}_X}(\mathcal{O}_X, \mathcal{O}_X),$$ where $\mathcal{D}_X:=\mathcal{D}iff(\mathcal{O}_X) \cong \mathscr{U}(\mathcal{O}_X,\mathcal{T}_X)$ is the sheaf of differential operators of $X$. 
	
	Thus, in this case we get a canonical isomorphism $$Cotor^\bullet_{\mathcal{J}_X}(\mathcal{O}_X, \mathcal{O}_X) \cong Ext^\bullet_{\mathcal{D}_X}(\mathcal{O}_X, \mathcal{O}_X),$$
	where $\mathcal{J}_X:=\mathscr{J}(\mathcal{O}_X, \mathcal{T}_X)$ is the usual sheaf of jets on $X$ \cite{DRV,BP}.
\end{Cor}

\section {Further Remarks} 
\subsection{Lie algebroids over Commutative ringed spaces} All the previous discussions extend naturally to the setting of Lie algebroids over Noetherian separated schemes or schemes of finite type over a field of characteristic zero (see \cite{MK, UB}). In this framework, the structure sheaf is the sheaf of regular functions on the scheme, and the Lie algebroid structure is defined in a way that aligns with algebraic geometry. Many foundational results, such as cohomological constructions and module-theoretic properties, remain valid in this general setting.
Furthermore, these results extend beyond schemes to the more general setting of Lie algebroids over commutative ringed spaces \cite{SA}. Ultimately, one may study Lie algebroids in the broader framework of ringed sites (see \cite{CV, CRV, DRV}). This generalization requires shifting from a geometric perspective to a more algebraic viewpoint, as ringed sites allow for a more flexible approach to defining Lie algebroids in formal geometry. In particular, this approach is well suited to contexts such as derived geometry and the theory of stacks,  where the classical scheme-theoretic framework may not be sufficient.
\subsection{Deformations and low-dimensional cohomology}
Cohomology plays a fundamental role in the deformation theory of algebraic and geometric structures. In the classical theory of Lie algebras and associative algebras, Chevalley--Eilenberg cohomology and Hochschild cohomology govern infinitesimal deformations, automorphisms, and obstruction theory. For Lie algebroids, the corresponding deformation theory is governed by Chevalley--Eilenberg--de Rham (hyper)cohomology and Hochschild (hyper)cohomology, with certain modifications arising from the absence of a genuine adjoint representation in general (see \cite{CM, CV, CRV, TS}). In this setting, suitable Gerstenhaber brackets play a central role in encoding the deformation theory and the associated obstruction calculus. Consequently, the low-dimensional cohomology groups admit the usual deformation-theoretic interpretation.
As in \cite{UB, PT, TS, CM, CV, CRV}, these interpretations remain valid for the sheaves of algebras associated with Lie algebroids over topological ringed spaces in the more general setting considered in this article (for e.g. relaxing the rank condition, base space is with singularities).

\subsection*{Acknowledgments}
The author acknowledges support from the C. S. I. R. fellowship grant during his Ph.D. at IIT Kanpur. 
	He wholeheartedly expresses his deepest gratitude to Dr. Ashis Mandal for his immense support and guidance. Additionally, he would like to thank Dr. Satyendra Kumar Mishra for his valuable comments and suggestions, which helped to improve the article.

\vspace{.4 cm}
{\bf Abhishek Sarkar},~
Department of Applied Sciences and Humanities, 
\\School of Computing,
MIT Art, Design and Technology University, Pune,
India.\\
e-mail:  abhisheksarkar49@gmail.com , abhishek.sarkar@mituniversity.edu.in

\end{document}